\documentclass[a4paper,english, 12pt]{article}
\usepackage[latin9]{inputenc}
\usepackage{amstext}
\usepackage{amsthm}
\usepackage{graphicx}
\usepackage{amssymb}
\usepackage{tikz}
\usepackage{subfig}
\usepackage{amsmath}

\usepackage{float}
\usepackage[top = 3.5cm, bottom = 4cm, left = 2cm , right = 2cm]{geometry}

\numberwithin{equation}{section}
\numberwithin{figure}{section}
\theoremstyle{plain}
\newtheorem{thm}{\protect\theoremname}
\newtheorem{lem}[thm]{Lemma}
\theoremstyle{definition}
\newtheorem{defn}[thm]{\protect\definitionname}

\makeatother

\usepackage{babel}
\providecommand{\definitionname}{Definition}
\providecommand{\theoremname}{Theorem}

\usepackage{fancyhdr}
\usepackage{atbegshi}
\begin{document}

\title{\vspace{-2cm} Generalised Sierpinski Triangles}
\author{ Kyle Steemson  - Christopher Williams }
\date{ Australian National University \today}
\maketitle
\begin{center} Keywords: Sierpinski, Pedal triangles, Fractal Tiling \end{center}
\begin{abstract}
The family of Generalised Sierpinski triangles consist of the classical Sierpinski triangle, the previously well investigated  Pedal triangle and two new triangular shaped fractal objects denoted by $\triangle FNN$ and $\triangle FFN$. 
All of the generalised Sierpinski triangles are defined in terms of iterated functions systems (IFS's) found by generalising the classic IFS used for the Sierpinski triangle. 
In this paper the new IFSs for the two new types of fractal triangles are defined, the dimensions of the triangles are analysed, and applications for pedagogical use and tiling theory discussed.
\end{abstract}

\section{\bf Introduction} \label{intro}
The Sierpinski triangle is one of the most well known fractals. It is an object which has zero area and infinite boundary. It was first discovered in 1915 by Waclaw Sierpinski \cite{Sierpinski} and has been thoroughly researched since. \\

The study of the Sierpinski triangle is interesting for pedagogical reasons as instructors often refer to it because of its simplicity and historical significance in fractal geometry. The Sierpinski triangle offers students a simple way to understand self-similarity and fixed points outside the realm of $\mathbb{R}$. Students can also learn the properties of a fractal attractor and its construction via the chaos game through Sierpinski triangles with the Sierpinski triangle offering a simple $\mathbb{R}^2$ example of this. Generalised Sierpinski triangles are interesting for a similar reason because they offer an extension to the classical Sierpinski triangle with fewer symmetries.\\ 
The family of generalised Sierpinski triangles is a set of four triangle shaped attractors found by generalising the iterated function system (IFS) of the Sierpinski triangle. An IFS and an attractor are defined in Section \ref{Section2}. 
The definition of a generalised Sierpinski triangle is given in Section \ref{fractrisec} as a geometric object (Definition \ref{eq: FracTriGeoDef}). Section \ref{proofof4} proves that the family of generalised Sierpinski triangles is made of exactly four types: $\triangle NNN$ (Sierpinski),$\triangle FNN$, $\triangle FFN$ and $\triangle FFF$ (Pedal). \\

The currently known Pedal triangle is a sub-triangle found within a larger triangle by a geometrical method as follows: for the triangle $\triangle ABC$ and the selection of any point $P$ within $\triangle ABC$, we define $X$ such that the lines $AC$ and $XP$ are perpendicular. Points $Y$ and $Z$ are similarly defined. This generates $\triangle XYZ$ which is a Pedal triangle. Figure 1.1 shows one particular case.

\begin{figure}[H]
\centering 
\begin{tikzpicture}[scale = 0.7]
	\draw (0,0)--(2.49, 4.1036442)--(6,0)--cycle;	
	\node[below,left] at (0,0){$A$};
	\node[below,right] at (6,0){$B$};
	\node[above] at (2.49, 4.1036442){$C$};
	\draw (1.125047168,1.854134085)--(2.49, 0)--(4.01026875,2.326253602)--cycle;
	\draw (1.125047168,1.854134085)--(2.49, 1.02591118);
	\draw (2.49, 0)--(2.49, 1.02591118);
	\draw (4.01026875,2.326253602)--(2.49, 1.02591118);
	\node[left] at (1.125047168,1.854134085){$X$};
	\node[below] at (2.49, 0){$Y$};
	\node[right] at (4.01026875,2.326253602){$Z$};
	\node[above] at (2.49, 1.02591118){$P$};
\end{tikzpicture}
\caption{For $\triangle ABC$ the point $P$ can be chosen to generate the Pedal triangle $\triangle XYZ$.}
\label{fig: PedalP}
\end{figure}
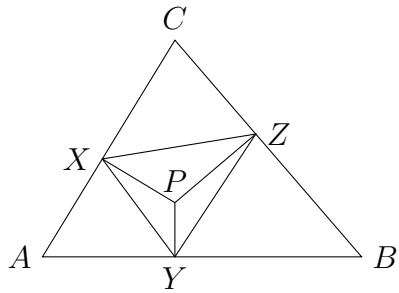

\begin{figure}[H]
	\centering
	\begin{tabular}{cccc}
		\subfloat[$\triangle NNN$ \label{fig:NNNBW}]
		{\includegraphics[width = 1.5in]{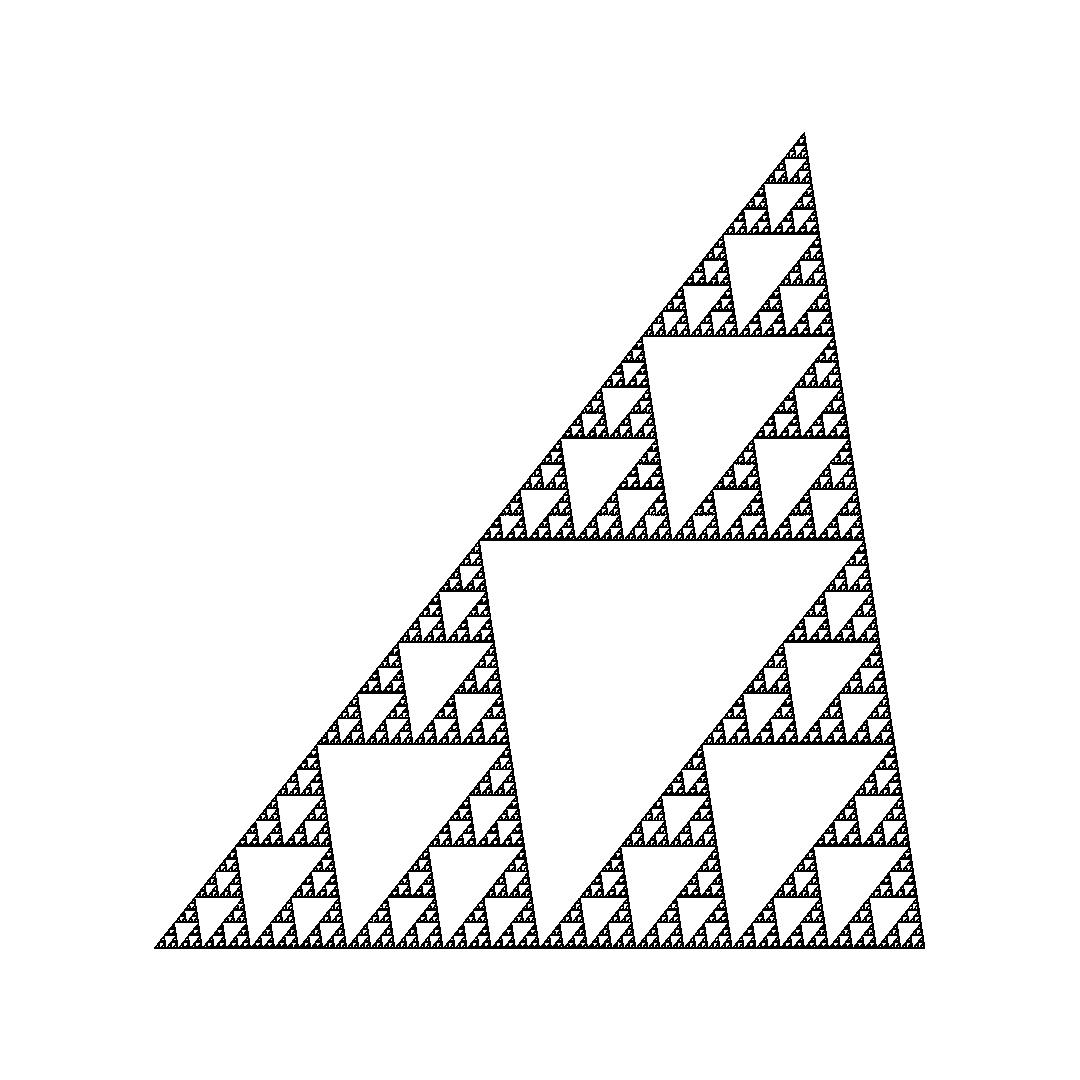}} &
		
		\subfloat[$\triangle FNN$ \label{fig:FNNBW}]
		{\includegraphics[width = 1.5in]{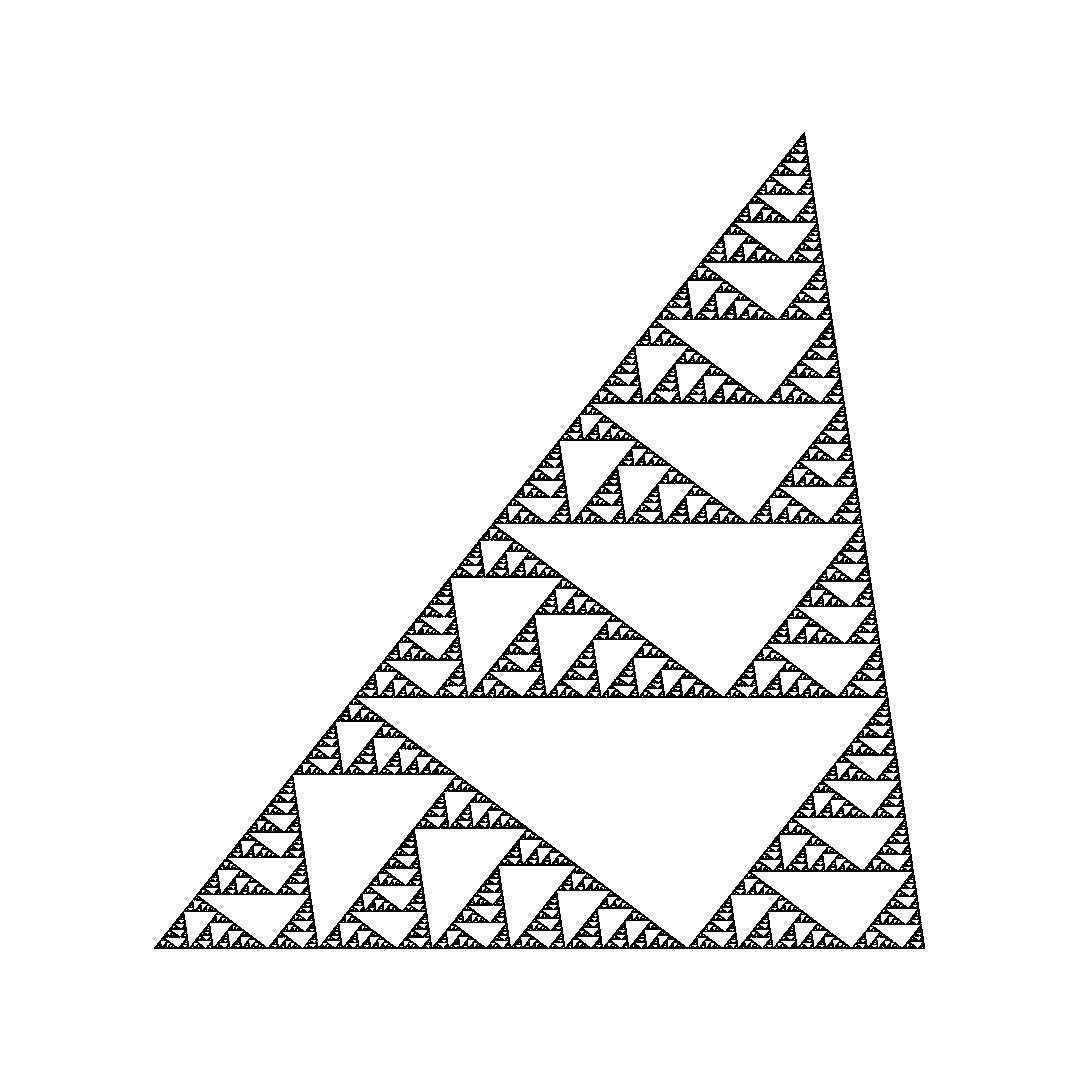}} &
		
		\subfloat[$\triangle FFN$ \label{fig:FFNBW}]
		{\includegraphics[width = 1.5in]{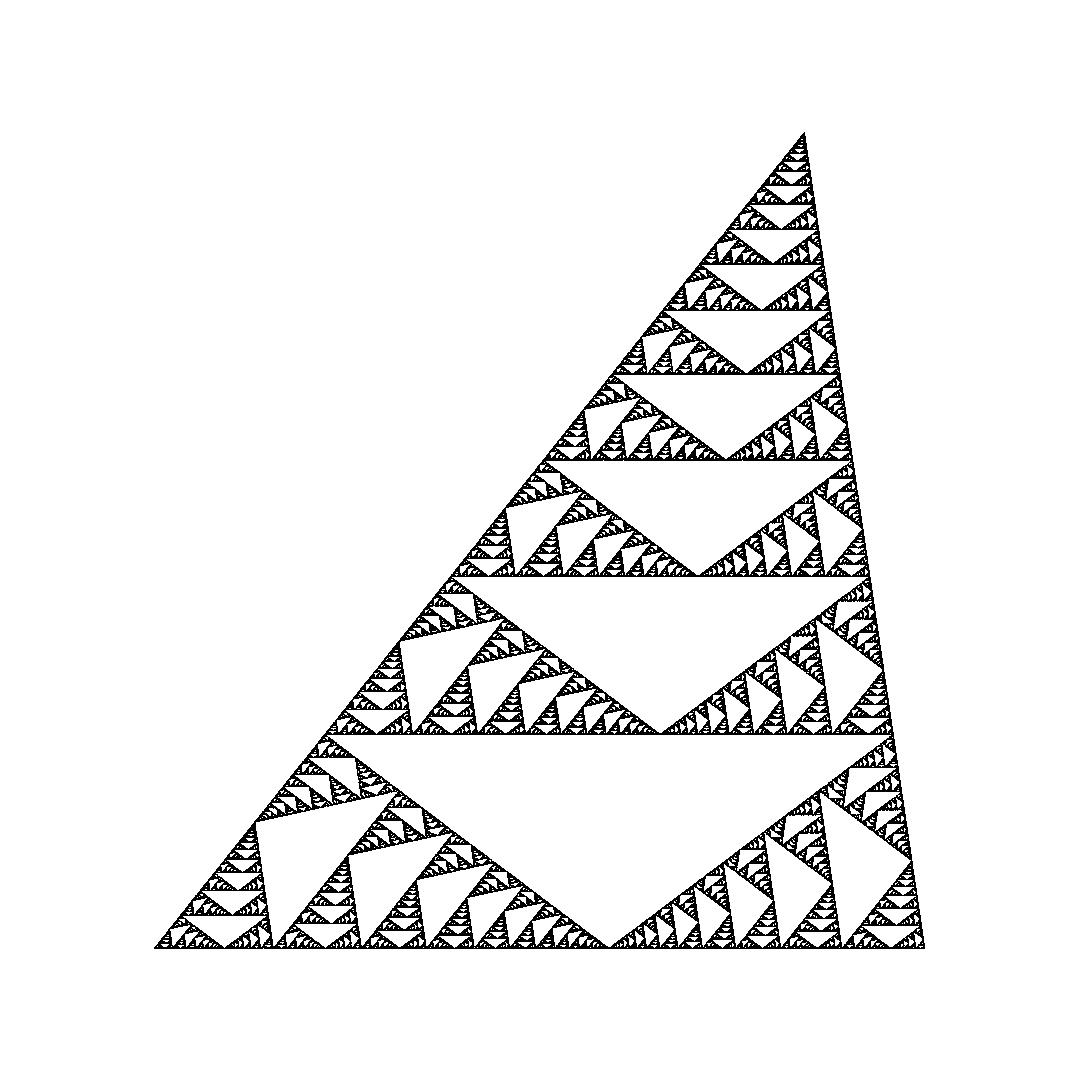}} &
		
		\subfloat[$\triangle FFF$ \label{fig:FFFBW}]
		{\includegraphics[width = 1.5in]{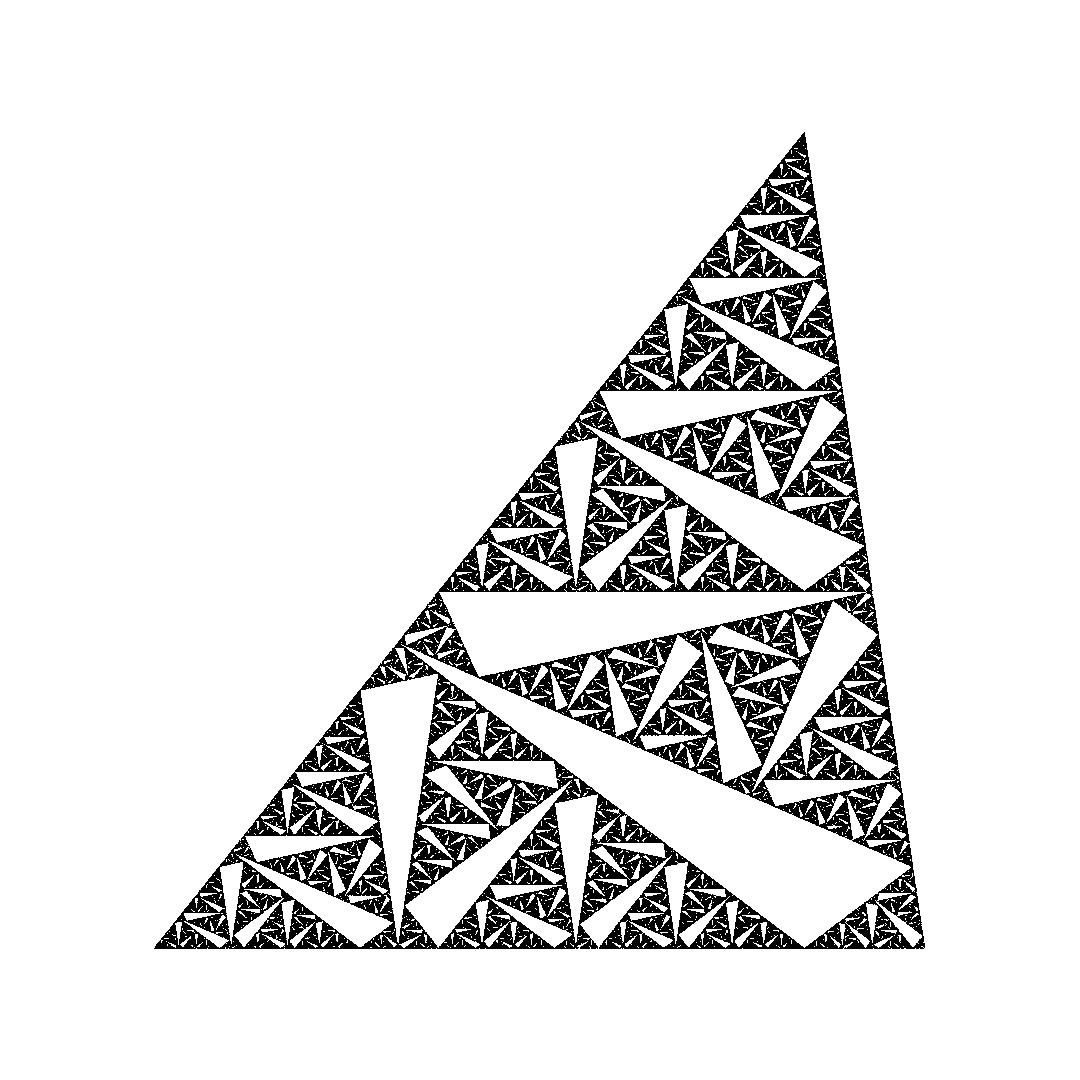}} 
	\end{tabular}
	\caption{The Four Generalised Sierpinski triangles corresponding to a particular fixed triangle}
	\label{fig: intropic}
\end{figure}

A fractal Pedal triangle is the attractor of an IFS such that $\triangle ABC$ is mapped into each of the three corner triangles by means of an affine transformation. For any $\triangle ABC$ there exists a Pedal triangle $\triangle XYZ$ such that the affine maps are similitudes. The fractal Pedal triangle and its dimension have been well researched in [8], [3], [4], and [5].\\
In addition to the Pedal triangle, two new generalisations of the Sierpinski triangle (denoted $\triangle FNN$ and $\triangle FFN$) are defined in Sections \ref{FNN Section} and \ref{FFN Section} respectively. It appears that these generalisations have not been considered in the literature before. This is potentially because the Pedal triangle is normally considered with respect to its geometrical construction and not in terms of an attractor of a generalised Sierpinski IFS.\\
It has been proven elsewhere that for any triangle there exists a unique Sierpinski triangle that uses the initial triangle as its convex hull. We prove in this paper that this is true for all generalisations of Sierpinski triangles . One example of this is seen in Figure \ref{fig: intropic}. The differences are most easily seen by analysing the size, shape and orientation of the white holes in the black fractal attractor.\\ 

\vspace{-0.5cm}

In Section \ref{fractrisec}, we give an IFS definition for each type of generalised Sierpinski triangle in terms of the triangle side lengths . In Section \ref{tridim} the fractal (Hausdorff) dimension of each of the generalised Sierpinski triangles is investigated. Fractal tilings are also briefly discussed. The tiling theory developed in \cite{BarnNew} has been used to create an example of a periodic generalised Sierpinski triangle tiling in Section \ref{fractalblowups}.   \\

\section{ \bf Basic Fractal Geometry Theory} \label{Section2} 
It is necessary to first build the fundamental mathematics for the family of generalised Sierpinski triangles. 
The following mathematics could be generalised to a complete metric space $(\mathbb{X},d)$,but for this paper it is sufficient to work in $(\mathbb{R}^2,d_{Euclid})$ which is a complete metric space. 
\begin{defn}
Iterated Function System (IFS)\\
If $f_{i}:\mathbb{R}^2 \rightarrow \mathbb{R}^2$, $i=1,2,...,N$,  are contractive functions, then $F:=(\mathbb{R}^2;f_{1},f_{2},....,f_{N})$ is called an $\mathit{iterated}$ $\mathit{function}$ $\mathit{system}$.\\
Contractive means that $d(f_{i}(x)-f_{i}(y))\leq\lambda_{i}d(x-y)$ $\forall$$x,y\in \mathbb{R}^2$, $\lambda_{i}\in[0,1)$, and $i=1, ..., N$. 
\end{defn}
Each of the $f_i$ are affine maps; that is, they are of the form $f_i(x)=L_i x + b_i$ where $L_i$ is a linear transformation given by a $2 \times 2$ real matrix and $x$ and $b_i$ are $2 \times 1$ vectors. If $L_i$ can be expressed as $L_i= \lambda_i O_i$ where $O_i$ is a $2 \times 2$ orthogonal matrix then $f_i$ is a similitude. Generalised Sierpinski triangles are defined using similitudes in their IFSs. 
\begin{defn}
Invertible IFS\\
If each of the maps $f\in F$ are homeomorphisms then $F$ is said to be an $\mathit{invertible}$ IFS and $F^{-1}$ is defined by:  $F^{-1}:=(\mathbb{R}^2;f_{1}^{-1},f_{2}^{-1},....,f_{N}^{-1})$.
\end{defn}
All affine maps in $\mathbb{R}^2$ with $\mathrm{det}(L) \neq 0$ are invertible. Hence, an IFS of affine maps that have non-zero determinant is invertible. Of interest for this paper are the Hausdorff subsets (denoted $\mathbb{H(X)}$), which are the non-empty compact sets in $(\mathbb{R}^2, d_{Euclid})$ or equivalently are the closed and bounded subsets of $\mathbb{R}^2$.\\
By slight abuse of notation, $F$ is used to denote how an IFS acts on a Hausdorff subset of $\mathbb{X}$.

\begin{defn}The Operation of an IFS \\
An IFS $F:\mathbb{H}(\mathbb{R}^2)\rightarrow \mathbb{H}(\mathbb{R}^2)$ acting on a point $A\in\mathbb{H}(\mathbb{R}^2)$ is defined as, \begin{equation} F(A)= \bigcup_{f_{i}\in F}f_{i}(A) = f_1(A) \cup \cdots \cup f_N(A) \end{equation}.
\end{defn}

The notation $F^k$ is introduced for repeated application of $F$, that is $F^k(A)=F( \cdots F(A) \cdots )$ ($k$ times) for some Hausdorff subset $A$ and $k \geq 1$.

\begin{defn} Attractor of an IFS\\
A non-empty set $A\in\mathbb{H}(\mathbb{R}^2)$ is the attractor of an IFS $F$ if both of the following conditions are satisfied:
\begin{enumerate}
\item  $F(A)=A$, and
\item  There exists an open set $U \subset \mathbb{R}^2$ such that for all $S\in\mathbb{H}(\mathbb{R}^2)$ with $S \subset U$, we have $A\subset U$ and $\lim_{k\rightarrow\infty}F^{k}(S)=A$ , where the limit is taken with respect to the Hausdorff metric $d_{\mathbb{H}}$ on $\mathbb{H}(\mathbb{R}^2)$. 
\end{enumerate}
\end{defn}

The Hausdorff  metric is a metric between sets in $\mathbb{R}^2$ defined for the points  $X, Y \in \mathbb{H}(\mathbb{R}^2)$ as, 
\begin{equation*}
d_{\mathbb{H}}(X,Y) = \max \bigg \{ \max_{x \in X} \min_{y \in Y} \big \{ d_{Euclid}(x,y) \big \}, \max_{y \in Y} \min_{x \in X} \big \{ d_{Euclid}(x,y) \big \} \bigg \}
\end{equation*}

This next theorem is an extension of Banach's contraction mapping theorem for the case when the function is mapping from $\mathbb{H}(\mathbb{R}^2)$ to $\mathbb{H}(\mathbb{R}^2)$.  

\begin{thm} Hutchinson's Theorem \cite{Hutch} \\
If an IFS $F$ is contractive on a complete metric space $\mathbb{{X}}$, then F has a unique attractor with the whole space $\mathbb{{X}}$ being the basin of attraction. 
\end{thm}
Therefore, in $(\mathbb{R}^2, d_{Euclid})$ with contractive IFS $F$ there is a unique attractor and all of $\mathbb{R}^2$ is the basin of attraction. For this paper, the unique attractor will be the generalised Sierpinski triangles. 

\begin{defn} The Open Set Condition (OSC) \\
The IFS $F = (\mathbb{R}^2, f_1, \cdots f_N)$ obeys the open set condition if, 
\begin{enumerate}
\item There exists a non-empty open set $\mathcal{O}$ such that $f_i(\mathcal{O}) \subset \mathcal{O}$ for all $i=1, \cdots, N$, and
\item $f_i(\mathcal{O}) \cap f_j(\mathcal{O}) = \emptyset$ for all $i,j \in \{1, \cdots, N \}$ with $i \neq j$.
\end{enumerate}
\label{OSC}
\end{defn}

\begin{thm}The Moran--Hutchinson Theorem \cite{Hutch} \\
If an attractor is constructed from an IFS, if it satisfies the OSC and if its maps are contractive similitudes of scaling factors $0<s_i<1$, then the fractal dimension is the unique positive solution $d$ to $\sum_{i=1}^{N} s_i^d = 1$
\label{hutmor}
\end{thm}
For the family of generalised Sierpinski triangles, the IFSs are comprised of similitudes and satisfy the OSC. Therefore the fractal dimension can be calculated using the Moran--Hutchinson theorem.

\section{Generalised Sierpinski Triangles} \label{fractrisec}
The aforementioned (see Section \ref{intro}) four types of generalised Sierpinski triangles are distinguished by the types of mappings they include. A flip (F) map refers to a contractive similitude that involves a flip along an axis, rotation, scaling and translation. Geometrically a flip map performed in $\mathbb{R}^2$ on a triangle will keep one corner of the triangle fixed and exchange the lengths of the two adjacent sides. A non-flip (N) map refers to a contractive similitude that only scales and translates. 

\begin{defn} \label{eq: FracTriGeoDef} Generalised Sierpinski Triangles \\
An attractor of an IFS is called a generalised Sierpinski triangle when it satisfies the following,
\begin{enumerate}
\item The IFS consists of exactly three continuous and contractive similitudes.
\item The fixed points for each map form the vertices of the triangle.
\item The attractor is non overlapping and just touching on the three points where the maps meet and each of these points are collinear with two different vertices.
\end{enumerate}
\end{defn}

The IFS having contractive maps (Property 1) implies that an attractor must always exist by the Hutchinson Theorem and additionally implies that similitudes give a calculable fractal dimension by the Moran--Hutchinson Theorem. This allows investigation into fractal tilings with the different generalised Sierpinski triangles. Property 2 and 3 geometrically describe the result of Property 1.  The vertices are required to be fixed points because otherwise the maps are not similitudes. The set of overlap between the maps must also be collinear with the vertices as the IFS is constructed from similitudes. \\

A proof of why $\triangle NNN, \triangle FNN, \triangle FFN$ and $\triangle FFF$ is an exhaustive list of the generalised Sierpinski triangles is provided in Section \ref{proofof4} following the explicit definition of each type of IFS. \\

\begin{figure}[h!]
\centering 
\begin{tikzpicture}[scale=0.5]
	\draw (0,0)--(2.4, 6.78)--(8,0)--cycle;	
	\node[below,left] at (0,0){$A (0,0)$};
	\node[below,right] at (8,0){$B (1,0)$};
	\node[above] at (2.4, 6.78){C $(C_x,C_y)$};
	\draw[<->] (-0.5,0.5)--(1.8, 6.78)
		node[pos=0.5,fill=white]{$b$};
	\draw[<->] (8.5,0.5)--(3.2, 6.78)
		node[pos=0.5,fill=white]{$a$};
	\draw[<->] (0, -1.3)--(8,-1.3)
		node[pos=0.5,fill=white]{$1$};
	\draw[dashed] (2.4,0) -- (2.4,6.78)
	node[pos=0.5,right]{$C_y$};
	\draw[<->] (0,-0.7)--(2.4,-0.7)
	node[pos=0.5,fill=white]{$C_x$};
\end{tikzpicture}
\caption{The zeroth order (convex hull) of all generalised Sierpinski triangles can be resized so that it is defined by $a$ and $b$, for this diagram $a=1.1$ and $b=0.9$.}
\label{fig: CxCytri}
\end{figure}
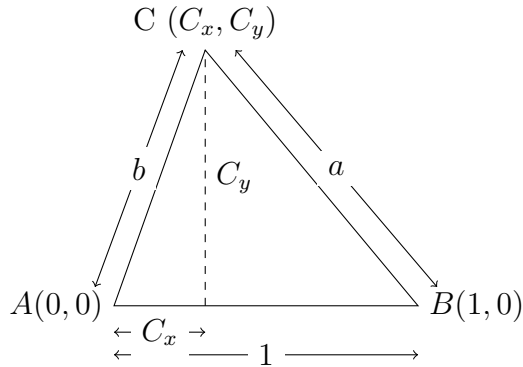

Without loss of generality, the fractal triangle $\triangle ABC$ can be rescaled, rotated and translated such that the bottom side be of length one and run from $(0,0)$ to $(1,0)$ in $\mathbb{R}^2$. Performing this linear transform to any fractal triangle reduces the number of parameters without altering any properties of the fractal. 
The exterior of any fractal triangle is now fully defined by the side lengths $a$ and $b$ as seen in Figure \ref{fig: CxCytri}, and this paper will show that the IFSs for each of the four types of generalised Sierpinski triangles can be purely parametrised by $a$ and $b$. 

Through elementary geometric calculation, the coordinates of the vertex $C(C_x, C_y)$ are computed to be:
\begin{equation}
C(C_x, C_y) =  \Big ( \frac{1}{2}( b^2 - a^2+1) , \frac{1}{2} \sqrt{4b^2 - ( b^2 - a^2+1)^2}  \Big )
\end{equation}

\newpage
\subsection{Sierpinski ($\bf \triangle NNN$) Triangles} 
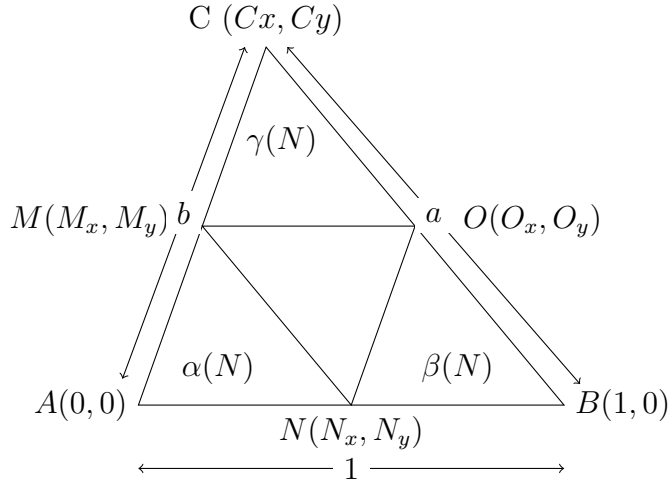
\begin{figure}[h!] 
\centering 
\begin{tikzpicture}[scale=0.7]
	\draw (0,0)--(2.4, 6.78)--(8,0)--cycle;	
	\node[below,left] at (0,0){$A (0,0)$};
	\node[below,right] at (8,0){$B (1,0)$};
	\node[above] at (2.4, 6.78){C $(Cx,Cy)$};
	\draw (1.2,3.39)--(4,0)--(5.2,3.39)--cycle;
	\node[left] at (0.8,3.5){$M(M_x, M_y)$};
	\node[below] at (4,0){$N(N_x, N_y)$};
	\node[right] at (5.9,3.5){$O(O_x,O_y)$};
	\node[below] at (1.5,1.2){$\alpha (N)$};
	\node[below] at (6,1.2){$\beta (N)$};
	\node[below] at (2.7,5.5){$\gamma (N)$};
	\draw[<->] (-0.3,0.5)--(2, 6.78)
		node[pos=0.5,fill=white]{$b$};
	\draw[<->] (8.3,0.4)--(2.8, 6.78)
		node[pos=0.5,fill=white]{$a$};
	\draw[<->] (0, -1.2)--(8,-1.2)
		node[pos=0.5,fill=white]{$1$};
\end{tikzpicture}
\caption{First Order Sierpinski ($\triangle NNN$) triangle with $a=1.1$ and $b=0.9$}
\label{fig: NNNtri}
\end{figure}

\hspace*{\fill} \\
The three non-flip similitudes are described by where the vertices of the large triangle are mapped, see Figure \ref{fig: NNNtri}.
\begin{equation}
\begin{aligned}[c]
f_A(0,0) &= (0,0) \\ f_B(0,0) &= (N_x,N_y) \\ f_C(0,0) &= (M_x,M_y)
\end{aligned} \qquad
\begin{aligned}[c]
f_A(1,0) &= (N_x, N_y) \\ f_B(1,0) &= (1,0) \\  f_C(1,0) &= (O_x,O_y)
\end{aligned} \qquad
\begin{aligned}[c]
f_A(C_x, C_y) &= (M_x, M_y) \\  f_B(C_x, C_y) &= (O_x, O_y) \\  f_C(C_x, C_y) &= (C_x, C_y) 
\end{aligned}
\end{equation}

The scaling ratios of maps $f_A, f_B, f_C$ are $\alpha, \beta, \gamma$ respectively. Three equations which describe how the scaling ratios are related are found from the three sides of the fractal triangle. 
\begin{equation} \label{ref: eq 5}
\begin{aligned}[c]
\alpha \cdot 1 + \beta \cdot 1 &= 1 \\ \alpha \cdot b + \gamma \cdot b &= b \\ \beta \cdot a + \gamma \cdot a &= a 
\end{aligned} \qquad \Rightarrow \qquad 
\begin{aligned}[c]
\alpha + \beta &= 1 \\ \alpha + \gamma &= 1 \\ \beta + \gamma &= 1
\end{aligned} 
\end{equation}
The solution to Equation \ref{ref: eq 5} is $\alpha = \beta = \gamma = \frac{1}{2}$. Therefore the scaling ratio for each map is always $\frac{1}{2}$ regardless of the side lengths $a$ or $b$  of the Sierpinski fractal triangle.

Note that the three points $N,M,O$ are the set of overlap of the IFS and themselves define a triangle-shaped hole in the fractal. For the Sierpinski case ($\triangle NNN$), it is important to realise that each side of the small centre triangle is parallel to a unique side of the large triangle. This property will be lost with the introduction of flip maps into the IFS.

The translation for each map in the IFS is given by where the origin is mapped. For $f_A$ the origin is the fixed point so it does not have a translation. For $f_B$ the origin is mapped to $(N_x, N_y) = (\alpha,0) = (\frac{1}{2},0)$. For $f_C$ the origin is mapped to $(M_x,M_y) = (\frac{\alpha b}{b}C_x, \frac{\alpha b}{b}C_y ) = (\frac{C_x}{2}, \frac{C_y}{2})$. Therefore the IFS for the classical Sierpinski triangle contains the following maps from $\mathbb{R}^2$ to $\mathbb{R}^2$, 

\begin{equation} f_A 
\begin{bmatrix}
    x  \\      y  \\
\end{bmatrix} =
\begin{bmatrix}
    \frac{1}{2}       & 0  \\      0   & \frac{1}{2}  
\end{bmatrix}  \begin{bmatrix}
    x  \\      y  \\
\end{bmatrix}  \end{equation}
\begin{equation}  f_B 
\begin{bmatrix}
    x  \\      y  \\
\end{bmatrix}  =
\begin{bmatrix}
    \frac{1}{2}       & 0  \\      0   & \frac{1}{2}  
\end{bmatrix} \begin{bmatrix}
    x  \\      y  \\
\end{bmatrix} +
\begin{bmatrix}
    \frac{1}{2}  \\     0  \\
\end{bmatrix} \end{equation}
\begin{equation} f_C 
\begin{bmatrix}
    x  \\     y  \\
\end{bmatrix} =
\begin{bmatrix}
    \frac{1}{2}       & 0  \\     0   & \frac{1}{2}  
\end{bmatrix} \begin{bmatrix}
    x  \\     y  \\
\end{bmatrix} +
\begin{bmatrix}
    \frac{C_x}{2}  \\     \frac{C_y}{2}  \\
\end{bmatrix} \end{equation}

The maps $f_A$ and $f_B$ are completely independent of the shape of the Sierpinski triangle; only the translation in the map $f_C$ depends on the top vertex $C$ (which in turn depends on the side lengths $a$ and $b$). Hence, $a$ and $b$ completely determine the IFS that has the Sierpinski triangle as the attractor.

\vspace{2mm}
\subsection{\bf Generalised Sierpinski Triangle: $\bf \triangle FNN$} \label{FNN Section} \hspace*{\fill} \\
It is assumed that the flip map of the $\triangle FNN$ IFS is $f_A$ with scaling ratio $\alpha$ and that the non-flip maps are $f_B$ and $f_C$ with scaling ratios $\beta$ and $\gamma$ respectively. 
The three functions are described by where the vertices of the large triangle are mapped. 
\begin{equation} \label{eq: FNNmaps}
\begin{aligned}[c]
f_A(0,0) &= (0,0) \\ f_B(0,0) &= (N_x,N_y) \\ f_C(0,0) &= (M_x,M_y)
\end{aligned} \qquad
\begin{aligned}[c]
f_A(1,0) &= (M_x, M_y) \\ f_B(1,0) &= (1,0) \\  f_C(1,0) &= (O_x,O_y)
\end{aligned} \qquad
\begin{aligned}[c]
f_A(C_x, C_y) &=  (N_x, N_y) \\  f_B(C_x, C_y) &= (O_x, O_y) \\  f_C(C_x, C_y) &= (C_x, C_y) 
\end{aligned}
\end{equation}

\begin{figure}[h!]
\centering 
\begin{tikzpicture}[scale=0.7]
	\draw (0,0)--(2.4, 6.78)--(8,0)--cycle;	
	\node[below,left] at (0,0){$A (0,0)$};
	\node[below,right] at (8,0){$B (1,0)$};
	\node[above] at (2.4, 6.78){C $(Cx,Cy)$};
	\draw (1.32,3.75)--(3.58,0)--(4.9,3.75)--cycle;
	\node[left] at (0.8,3.8){$M(M_x, M_y)$};
	\node[below] at (3.58,0){$N(N_x, N_y)$};
	\node[right] at (5.7,3.8){$O(O_x,O_y)$};
	\node[below] at (1.5,1.2){$\alpha (F)$};
	\node[below] at (6,1.2){$\beta (N)$};
	\node[below] at (2.7,5.5){$\gamma (N)$};
	\draw[<->] (-0.3,0.5)--(2, 6.78)
		node[pos=0.5,fill=white]{$b$};
	\draw[<->] (8.3,0.4)--(2.8, 6.78)
		node[pos=0.5,fill=white]{$a$};
	\draw[<->] (0, -1.2)--(8,-1.2)
		node[pos=0.5,fill=white]{$1$};
\end{tikzpicture}
\caption{First Order FNN fractal triangle with $a=1.1$ and $b=0.9$}
\label{fig: FNNtri}
\end{figure}

Each side of the fractal triangle in Figure \ref{fig: FNNtri} gives an equation which relates the scaling ratios of the IFS maps. 
\begin{equation} \label{FNNscaling} \begin{aligned}
\alpha.b + \beta.1 &= 1 \\ \alpha.1 + \gamma.b &= b \\ \beta.a + \gamma.a &= a
\end{aligned} \qquad \Rightarrow \qquad 
\begin{bmatrix}
    \alpha \\ \beta \\ \gamma
\end{bmatrix} =  \begin{bmatrix}
   \frac{b}{b^2+1} \\ \frac{1}{b^2+1} \\ \frac{b^2}{b^2+1}
\end{bmatrix} \end{equation}
Note that these maps only depend on the side length $b$. 
Now the translations for each map can be found. $f_A$ has no translation. $f_B$ maps the origin to $(N_x,N_y) = (\alpha b, 0) = (1-\beta,0)$ so this is the translation. $f_C$ maps the origin to $(M_x, M_y) = (\frac{\alpha}{b}C_x, \frac{\alpha}{b}C_y) = ( (1-\gamma)C_x, (1-\gamma)C_y) )$, where relations from Eq \ref{FNNscaling} have been used. 
Therefore the non-flip maps $f_B$ and $f_C$ are described by: 

\begin{equation} f_B \begin{bmatrix}
    x  \\      y  \\
\end{bmatrix}  = \begin{bmatrix}
    \beta       & 0  \\      0   & \beta  
\end{bmatrix} \begin{bmatrix}
    x  \\      y  \\
\end{bmatrix} + \begin{bmatrix}
   1-\beta  \\     0 
\end{bmatrix} \end{equation}
\begin{equation} f_C \begin{bmatrix}
    x  \\     y  \\
\end{bmatrix} = \begin{bmatrix}
    \gamma       & 0  \\     0   & \gamma  
\end{bmatrix} \begin{bmatrix}
    x  \\     y  \\
\end{bmatrix} + \begin{bmatrix}
    (1-\gamma)C_x  \\     (1-\gamma)C_y  \\
\end{bmatrix} \end{equation}

The flip map $f_A$ could be calculated using the Fundamental Theorem of Affine Geometry \cite{BarnTile} but instead a geometric argument is used. The map is comprised of a flip along the $x$-axis, a clockwise rotation by angle $A$ and scaling by $\alpha$; therefore the map must take the form of
\begin{equation} \begin{aligned} f_A \begin{bmatrix}
    x  \\      y  \\
\end{bmatrix} & = \begin{bmatrix}
    \alpha       & 0  \\      0   & \alpha  
\end{bmatrix}\begin{bmatrix}
    \cos(-A)       & -\sin(-A)  \\      \sin(-A)   & \cos(-A)
\end{bmatrix}\begin{bmatrix}
    1       & 0  \\      0   & -1  
\end{bmatrix} \begin{bmatrix}
    x  \\      y  \\
\end{bmatrix} \\  &= \begin{bmatrix}
   \alpha \cos(A)       & -\alpha \sin(A)  \\     -\alpha \sin(A)   & -\alpha \cos(A)
\end{bmatrix}\begin{bmatrix}
    x  \\      y  \\
\end{bmatrix} \\ &= \begin{bmatrix}
   q   &  p  \\     p  & -q
\end{bmatrix}\begin{bmatrix}
    x  \\      y  \\
\end{bmatrix} \end{aligned} \end{equation}
For some $p$ and $q$ in $\mathbb{R}$. These values can be found since $f_A(1,0) = (M_x, M_y) = (\frac{\alpha}{b}C_x, \frac{\alpha}{b}C_y)$, so:
\begin{equation} \begin{bmatrix}
    \frac{\alpha}{b}C_x  \\     \frac{\alpha}{b}C_y \\
\end{bmatrix} =  f_A \begin{bmatrix}
    1  \\      0  \\
\end{bmatrix} = \begin{bmatrix}
   q   &  p  \\     p  & -q
\end{bmatrix}\begin{bmatrix}
    1  \\      0  \\
\end{bmatrix} = \begin{bmatrix}
    q  \\     p  \\
\end{bmatrix}\end{equation}
Therefore, $q = \frac{\alpha}{b}C_x$ and $p = \frac{\alpha}{b}C_y$. These values can be checked using the mapping of the top vertex, $f_A(C_x, C_y) = (M_x, M_y) = (\alpha b, 0)$ and noting that by Figure \ref{fig: CxCytri},  $C_x^2 + C_y^2=b^2$.
\begin{equation} \begin{aligned}  f_A  \begin{bmatrix}
    C_x  \\      C_y  \\
\end{bmatrix} & = \begin{bmatrix}
   \frac{\alpha}{b}C_x   &  \frac{\alpha}{b}C_y  \\    \frac{\alpha}{b}C_y  & -\frac{\alpha}{b}C_x
\end{bmatrix} \begin{bmatrix}
    C_x  \\      C_y  \\
\end{bmatrix} = \begin{bmatrix}
    \frac{\alpha}{b}C_x^2 + \frac{\alpha}{b}C_y^2  \\      \frac{\alpha}{b}C_xC_y - \frac{\alpha}{b}C_xC_y  \\
\end{bmatrix}\\  & = \begin{bmatrix}
    \frac{\alpha}{b} \Big ( C_x^2 + C_y^2 \Big)  \\     0  \\
\end{bmatrix}= \begin{bmatrix}
    \frac{\alpha}{b} b^2  \\     0  \\
\end{bmatrix}= \begin{bmatrix}
    \alpha b  \\     0  \\
\end{bmatrix}\end{aligned}\end{equation}
This confirms the values for $q$ and $p$; therefore the flip map $f_A$ is: 
\begin{equation} f_A \begin{bmatrix}
    x  \\      y  \\
\end{bmatrix}  = \begin{bmatrix}
    \frac{\alpha}{b}C_x       & \frac{\alpha}{b}C_y  \\      \frac{\alpha}{b}C_y   & -\frac{\alpha}{b}C_x
\end{bmatrix} \begin{bmatrix}
    x  \\      y  \\
\end{bmatrix} \end{equation}
For any given triangle with side lengths $a$ and $b$, the mappings for an IFS have been determined, so the attractor is a FNN fractal triangle. 

\vspace{2mm}
\subsection{\bf Generalised Sierpinski Triangle: $\bf \triangle FFN$} \label{FFN Section} \hspace*{\fill} \\
For $\triangle FFN$ generalised Sierpinski triangles, the flip maps are $f_A$ and $f_B$ with scaling factors $\alpha$ and $\beta$ respectively, while the non-flip map is $f_C$ with scaling factor $\gamma$. 

The three maps of the FFN triangle are described by where the vertices of the exterior triangle are mapped. 
\begin{equation} \label{eq: FFNmaps}
\begin{aligned}[c]
f_A(0,0) &= (0,0) \\ f_B(0,0) &= (O_x,O_y) \\ f_C(0,0) &= (M_x,M_y)
\end{aligned} \qquad
\begin{aligned}[c]
f_A(1,0) &= (M_x, M_y) \\ f_B(1,0) &= (1,0) \\  f_C(1,0) &= (O_x,O_y)
\end{aligned} \qquad
\begin{aligned}[c]
f_A(C_x, C_y) &=  (N_x, N_y) \\  f_B(C_x, C_y) &= (N_x, N_y) \\  f_C(C_x, C_y) &= (C_x, C_y) 
\end{aligned}
\end{equation}

\begin{figure}[h!]
\centering 
\begin{tikzpicture}[scale=0.7]
	\draw (0,0)--(2.4, 6.78)--(8,0)--cycle;	
	\node[below,left] at (0,0){$A (0,0)$};
	\node[below,right] at (8,0){$B (1,0)$};
	\node[above] at (2.4, 6.78){C $(Cx,Cy)$};
	\draw (1.18,3.36)--(3.2,0)--(5.23,3.36)--cycle;
	\node[left] at (0.7,3.4){$M(M_x, M_y)$};
	\node[below] at (3.2,0){$N(N_x, N_y)$};
	\node[right] at (5.8,3.4){$O(O_x,O_y)$};
	\node[below] at (1.5,1.2){$\alpha (F)$};
	\node[below] at (6,1.2){$\beta (F)$};
	\node[below] at (2.7,5.5){$\gamma (N)$};
	\draw[<->] (-0.3,0.5)--(2, 6.78)
		node[pos=0.5,fill=white]{$b$};
	\draw[<->] (8.3,0.4)--(2.8, 6.78)
		node[pos=0.5,fill=white]{$a$};
	\draw[<->] (0, -1.2)--(8,-1.2)
		node[pos=0.5,fill=white]{$1$};
\end{tikzpicture}
\caption{First Order FFN fractal triangle with $a=1.1$ and $b=0.9$}
\label{fig: FFNtri}
\end{figure}
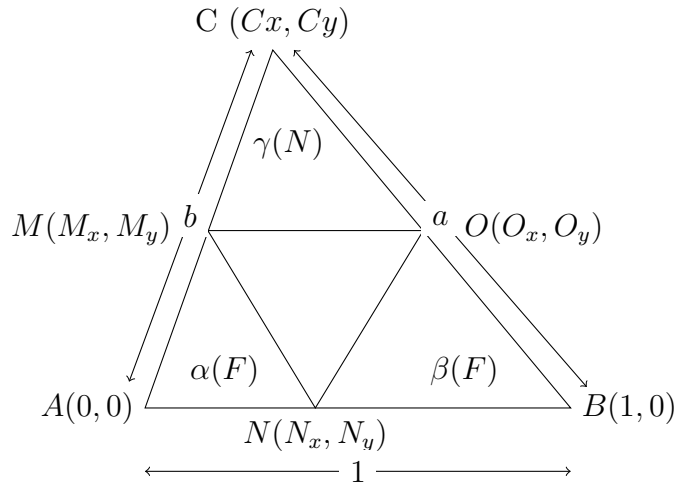

Each side of the fractal triangle gives an equation which relates the scaling ratios together. These equations are, 
\begin{equation} \begin{aligned}
\alpha \cdot b + \beta \cdot a &= 1 \\ \alpha \cdot 1 + \gamma \cdot b &= b \\ \beta \cdot 1 + \gamma \cdot a &= a
\end{aligned} \qquad \Rightarrow \qquad 
\begin{bmatrix}
    \alpha \\ \beta \\ \gamma
\end{bmatrix} =  \begin{bmatrix}
   \frac{b}{a^2+b^2} \\ \frac{a}{a^2+b^2} \\ \frac{a^2+b^2-1}{a^2+b^2}
\end{bmatrix} \end{equation}
In the case of a FFN fractal triangle the scaling ratios depend on both side lengths $a$ and $b$.
Now the translations for each map can be found. $f_A$ has no translation. $f_B$ maps the origin to $(O_x, O_y) = (1-\frac{\beta}{a}(1-C_x), \frac{\beta}{a}C_y))$ so this is the translation for $f_B$. $f_C$ maps the origin to $(M_x, M_y) = (\frac{\alpha}{b}C_x, \frac{\alpha}{b}C_y) = ( (1-\gamma)C_x, (1-\gamma)C_y) )$. 
The non-flip map $f_C$ is simply: 
\begin{equation} f_C \begin{bmatrix}
    x  \\     y  \\
\end{bmatrix} = \begin{bmatrix}
    \gamma       & 0  \\     0   & \gamma  
\end{bmatrix} \begin{bmatrix}
    x  \\     y  \\
\end{bmatrix} + \begin{bmatrix}
    (1-\gamma)C_x  \\     (1-\gamma)C_y  \\
\end{bmatrix} \end{equation}

From the FNN generalised Sierpinski triangles, we recall that flip maps take the form of: 
\begin{equation} f_{flip} \begin{bmatrix}
    x  \\      y  \\
\end{bmatrix} = \begin{bmatrix}
   q   &  p  \\     p  & -q
\end{bmatrix}\begin{bmatrix}
    x  \\      y  \\
\end{bmatrix}  + \begin{bmatrix}
    T_x  \\     T_y  \\
\end{bmatrix}\end{equation}
where $(T_x, T_y)$ is the translation of the map. 

Immediately, it is known that $f_A$ is given by,
\begin{equation} f_A \begin{bmatrix}
    x  \\      y  \\
\end{bmatrix} = \begin{bmatrix}
   \frac{\alpha}{b}C_x   &  \frac{\alpha}{b}C_y  \\     \frac{\alpha}{b}C_y  & -\frac{\alpha}{b}C_x
\end{bmatrix}\begin{bmatrix}
    \frac{\alpha}{b}C_x  \\     \frac{\alpha}{b}C_x  \\
\end{bmatrix} \end{equation}
This is the case because Eq. \ref{eq: FFNmaps} is identical to Eq. \ref{eq: FNNmaps} for the $f_A$ map and therefore the mapping that has already been calculated for FNN generalised Sierpinski triangles is correct. 

The $f_B$ mapping can be determined by noting that the origin is mapped to the point $(O_x, O_y)$, so the translation is $(T_x, T_y) = (1-\frac{\beta}{a}(1-C_x), \frac{\beta}{a}C_y))$. Now the values of $q$ and $p$ can be determined using the mapping $f_B(1,0) = (1,0)$. 
\begin{equation}  \begin{bmatrix}
    1 \\ 0 \\ 
\end{bmatrix} = f_B \begin{bmatrix}
    1 \\ 0 \\ 
\end{bmatrix} = \begin{bmatrix}
   q   &  p  \\     p  & -q
\end{bmatrix}\begin{bmatrix}
    1  \\      0  \\
\end{bmatrix}  + \begin{bmatrix}
    1-\frac{\beta}{a}(1-C_x) \\     \frac{\beta}{a}C_y  \\
\end{bmatrix} = \begin{bmatrix}
   q+1-\frac{\beta}{a}(1-C_x) \\     p+\frac{\beta}{a}C_y  \\
\end{bmatrix}\end{equation}
\begin{equation} \Rightarrow \begin{bmatrix}
  q  \\  p  \\
\end{bmatrix} = \begin{bmatrix}
  \frac{\beta}{a}(1-C_x)  \\  -\frac{\beta}{a}C_y  \\  
\end{bmatrix}  \end{equation}
These values for $q$ and $p$ can be confirmed using the mapping $f_B(C_x,C_y) = (M_x,M_y) = (1- \beta a, 0)$. 
The map $f_B$ is given by, 
\begin{equation} f_B \begin{bmatrix}
    x \\ y \\ 
\end{bmatrix} = \begin{bmatrix}
   \frac{\beta}{a}(1-C_x)   & -\frac{\beta}{a}C_y  \\     -\frac{\beta}{a}C_y  & -\frac{\beta}{a}(1-C_x)
\end{bmatrix}\begin{bmatrix}
    x  \\      y  \\
\end{bmatrix}  + \begin{bmatrix}
    1-\frac{\beta}{a}(1-C_x) \\     \frac{\beta}{a}C_y  \\
\end{bmatrix} \end{equation}
Hence, $f_A, f_B, f_C$ are all completely and explicitly determined by $a$ and $b$. 

\vspace{2mm}
\subsection{\bf Pedal ($\bf \triangle FFF$) Triangles} \label{FFF Section} \hspace*{\fill} \\
As before, the three flip maps of the Pedal (FFF) triangle are described by where the vertices of the exterior triangle are mapped. 
\begin{equation} \label{eq: FFFmaps}
\begin{aligned}[c]
f_A(0,0) &= (0,0) \\ f_B(0,0) &= (O_x,O_y) \\ f_C(0,0) &= (O_x, O_y)
\end{aligned} \qquad
\begin{aligned}[c]
f_A(1,0) &= (M_x, M_y) \\ f_B(1,0) &= (1,0) \\  f_C(1,0) &= (M_x,M_y)
\end{aligned} \qquad
\begin{aligned}[c]
f_A(C_x, C_y) &=  (N_x, N_y) \\  f_B(C_x, C_y) &= (N_x, N_y) \\  f_C(C_x, C_y) &= (C_x, C_y) 
\end{aligned}
\end{equation}

\begin{figure}[h!]
\centering 
\begin{tikzpicture}[scale=0.7]
	\draw (0,0)--(2.4, 6.78)--(8,0)--cycle;	
	\node[below,left] at (0,0){$A (0,0)$};
	\node[below,right] at (8,0){$B (1,0)$};
	\node[above] at (2.4, 6.78){C $(Cx,Cy)$};
	\draw (0.89,2.5)--(2.4,0)--(4.76,3.93)--cycle;
	\node[left] at (0.5,2.6){$M(M_x, M_y)$};
	\node[below] at (2.4,0){$N(N_x, N_y)$};
	\node[right] at (5.2,4.2){$O(O_x,O_y)$};
	\node[below] at (1.1,1.0){$\alpha (F)$};
	\node[below] at (5.5,1.5){$\beta (F)$};
	\node[below] at (2.7,5.2){$\gamma (F)$};
	\draw[<->] (-0.3,0.5)--(2, 6.78)
		node[pos=0.5,fill=white]{$b$};
	\draw[<->] (8.3,0.4)--(2.8, 6.78)
		node[pos=0.5,fill=white]{$a$};
	\draw[<->] (0, -0.9)--(8,-0.9)
		node[pos=0.5,fill=white]{$1$};
\end{tikzpicture}
\caption{First Order Pedal (FFF) fractal triangle with $a=1.1$ and $b=0.9$}
\label{fig: FFFtri}
\end{figure}
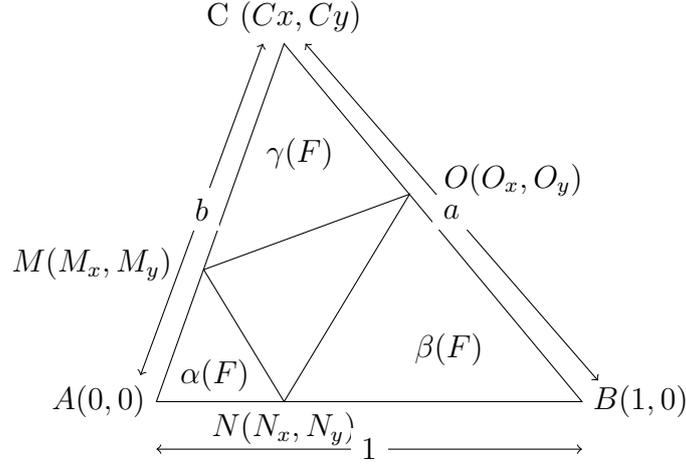

Each side of the fractal triangle gives an equation which relates the scaling ratios together. These equations are, 
\begin{equation*} \begin{aligned}
\alpha \cdot b + \beta \cdot a &= 1 \\ \alpha \cdot 1 + \gamma \cdot a &= b \\ \beta \cdot 1 + \gamma \cdot b &= a
\end{aligned} \qquad  \Rightarrow \qquad 
\begin{bmatrix}
    \alpha \\ \beta \\ \gamma
\end{bmatrix} =  \begin{bmatrix}
   \frac{1}{2b}(-a^2+b^2+1) \\ \frac{1}{2a}(a^2-b^2+1) \\ \frac{1}{2ab}(a^2+b^2-1)
\end{bmatrix} \end{equation*}

The translations for each map is where the origin is mapped. For $f_A$ the origin is a fixed point. For $f_B$ and $f_C$ the origin is mapped to $(O_x, O_y) = (1-\frac{\beta}{a}(C_x - 1), \frac{\beta}{a}C_y)$. The map for $f_A$ has already been determined in FNN and the map for $f_B$ has already been determined in FFN by noting the similarities between Eq \ref{eq: FFFmaps} with Eq \ref{eq: FFNmaps} and Eq \ref{eq: FNNmaps}. Since it is a flip map, $f_C$ must take the following form: 
\begin{equation} f_{flip} \begin{bmatrix}
    x  \\      y  \\
\end{bmatrix} = \begin{bmatrix}
   q   &  p  \\     p  & -q
\end{bmatrix}\begin{bmatrix}
    x  \\      y  \\
\end{bmatrix}  + \begin{bmatrix}
    T_x  \\     T_y  \\
\end{bmatrix}\end{equation}
where $(T_x, T_y)$ is the translation of the map so $(T_x, T_y)=(1-\frac{\beta}{a}(C_x - 1), \frac{\beta}{a}C_y)$. The point $(1,0)$ is mapped to $(M_x, M_y) = (\frac{\alpha}{b}C_x, \frac{\alpha}{b}C_y)$ so the values of $q$ and $p$ can be determined. 

\begin{equation}  \begin{bmatrix}
    \frac{\alpha}{b}C_x \\ \frac{\alpha}{b}C_y \\ 
\end{bmatrix} = f_C \begin{bmatrix}
    1 \\ 0 \\ 
\end{bmatrix} = \begin{bmatrix}
   q   &  p  \\     p  & -q
\end{bmatrix}\begin{bmatrix}
    1  \\      0  \\
\end{bmatrix}  + \begin{bmatrix}
    1-\frac{\beta}{a}(1-C_x) \\     \frac{\beta}{a}C_y  \\
\end{bmatrix} = \begin{bmatrix}
   q+1-\frac{\beta}{a}(1-C_x) \\     p+\frac{\beta}{a}C_y  \\
\end{bmatrix}\end{equation}
\begin{equation} \Rightarrow \begin{bmatrix}
  q  \\  p  \\
\end{bmatrix} = \begin{bmatrix}
  (\frac{\alpha}{b} - \frac{\beta}{a})C_x + (\frac{\beta}{a}-1)  \\  (\frac{\alpha}{b} - \frac{\beta}{a})C_y  \\  
\end{bmatrix}  \end{equation}
Therefore, given side lengths $a$ and $b$, the IFS for a fractal Pedal (FFF) triangle is: 

\begin{equation} f_A \begin{bmatrix}
    x  \\      y  \\
\end{bmatrix} = \begin{bmatrix}
   \frac{\alpha}{b}C_x   &  \frac{\alpha}{b}C_y  \\     \frac{\alpha}{b}C_y  & -\frac{\alpha}{b}C_x
\end{bmatrix}\begin{bmatrix}
    \frac{\alpha}{b}C_x  \\     \frac{\alpha}{b}C_x  \\
\end{bmatrix} \end{equation}
\begin{equation} f_B \begin{bmatrix}
    x \\ y \\ 
\end{bmatrix} = \begin{bmatrix}
   \frac{\beta}{a}(1-C_x)   & -\frac{\beta}{a}C_y  \\     -\frac{\beta}{a}C_y  & -\frac{\beta}{a}(1-C_x)
\end{bmatrix}\begin{bmatrix}
    x  \\      y  \\
\end{bmatrix}  + \begin{bmatrix}
    1-\frac{\beta}{a}(1-C_x) \\     \frac{\beta}{a}C_y  \\
\end{bmatrix} \end{equation}
\begin{equation} f_C \begin{bmatrix}
    x  \\      y  \\
\end{bmatrix} = \begin{bmatrix}
   (\frac{\alpha}{b} - \frac{\beta}{a})C_x + (\frac{\beta}{a}-1)   &  (\frac{\alpha}{b} - \frac{\beta}{a})C_y  \\     (\frac{\alpha}{b} - \frac{\beta}{a})C_y  & -(\frac{\alpha}{b} - \frac{\beta}{a})C_x - (\frac{\beta}{a}-1) \\
\end{bmatrix} \begin{bmatrix}
    x  \\      y  \\
\end{bmatrix}  + \begin{bmatrix}
    1-\frac{\beta}{a}(1-C_x)  \\     \frac{\beta}{a}C_y  \\
\end{bmatrix} \end{equation}

\vspace{2mm}
\subsection{Members of the Generalised Sierpinski Triangle Family} \label{proofof4} 
\begin{figure}[!htb] 
\centering 
\begin{tikzpicture}[scale=0.5]
	\draw[dotted] (0,0)--(2.4, 6.78)--(8,0)--cycle;	
	\node[below,left] at (0,0){$A(0,0)$};
	\node[below,right] at (8,0){$B (1,0)$};
	\node[above] at (2.4, 6.78){$C$};
	\draw (1.6,3.6)--(4,-0.6)--(4.9,3.2)--cycle;
	\draw(0,0)--(1.6,3.6)--(2.4,6.78)--(4.9, 3.2)--(8,0)--(4,-0.6)--cycle;
	\node[left] at (1.6,3.6){$M$};
	\node[below] at (4,-0.6){$N$};
	\node[right] at (4.9,3.2){$O$};
\end{tikzpicture}
\caption{Twisted Sierpinski Triangle - Proof of Collinear}
\label{fig: twisted Sierpinski}
\end{figure}
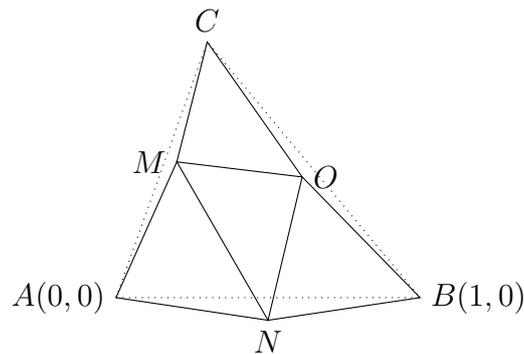

\noindent The claim that $\triangle NNN, \triangle FNN, \triangle FFN, \triangle FFF$ is the complete set of generalised Sierpinski triangles will be shown to be true as a corollary of the following theorem. \\
For the generalised Sierpinski triangles, Definition \ref{eq: FracTriGeoDef} requires that the set of overlap is collinear with the vertices. This will be proven to be a result of the generalised Sierpinski triangles being constructed from similitudes and will be done with respect to Figure \ref{fig: twisted Sierpinski}. \\

\begin{thm} \label{Kyle} For generalised Sierpinski Triangles, the set of overlap must be collinear with the vertices. \end{thm} 
\textit{Proof}: Given the fixed points $A,B,C$ of each of the IFS maps, there can be a triangle drawn that connects the vertices (the dashed semi-transparent shape in Figure \ref{fig: twisted Sierpinski}). There are three points of overlap in the Sierpinski triangle and this must also be the case in generalised Sierpinski triangles. The set of overlap is the points $N,M,O$ such that $f_A(A) = A, f_A( \{B,C\}) = \{N,M\}$ and similarly for $f_B$ and $f_C$, where $f_A, f_B, f_C$ are continuous contractive similitudes. 

In order to gain a contradiction, let us assume that $N$ is a point in the set of overlap and that $N$ is not collinear with the vertices $AB$. Suppose that $N$ sits underneath the $x$-axis. It is known by Lemma \ref{K} that angles are preserved by similitudes and that generalised Sierpinski triangles are required to be similitudes. Therefore the $\angle CAB$ is equal to $\angle MAN$ by invariance under $f_A$, and since $N$ is below the $x$-axis the point $M$ must be inside the interior of $\triangle ABC$. By an identical argument regarding $f_B$, the point $O$ must also be in the interior of $\triangle ABC$. Therefore $f_C$ must map $\angle ACB$ to $\angle MCO$ with both $M$ and $O$ being inside $\triangle ABC$. It is clear from Figure \ref{fig: twisted Sierpinski} that these two angles cannot be equal unless $\angle ACM$ and $\angle BCO$ both equal zero. This implies that the points $M$ and $O$ must in fact be collinear with $AC$ and $BC$ respectively. Thus $N$ is also collinear with $AB$ which contradicts the initial assumption. Therefore the set of overlap must be collinear with the fixed points as stated in Definition \ref{eq: FracTriGeoDef}. An identical argument would hold if $N$ was supposed to be above the x-axis instead of below. $\square$ \\

\newpage
\begin{lem} \label{K} Interior angles of a triangle are invariant under transformation by a similitude. \end{lem}
\textit{Proof:} The transformation of a similitude with scaling factor $s$ acting on a line in $\mathbb{R}^2$ space produces a line of which length has been scaled by a factor of $s$. This implies that the action of the similitude on a triangle in $\mathbb{R}^2$ produces a triangle with each side length scaled by a factor of $s$. Through a simple application of the cosine rule it is clear that the interior angles of a triangle are unchanged by a similitude. \\ 

An immediate consequence of Theorem \ref{Kyle} is that the only members of the family of generalised Sierpinski triangles is $\triangle NNN$, $\triangle FNN$, $\triangle FFN$ and $\triangle FFF$. This is because the only two types of similitude that conserve the fixed point and have the set of overlap collinear with the vertices is the non-flip (N) map and the flip (F) map. Therefore the possible combinations of these maps are the only members of the family.

\vspace{6mm}
\section{\bf Examples of Generalised Sierpinski Triangles }
\subsection{Simple construction} 
\begin{figure}[H]
	\centering
	\begin{tabular}{cccc}
		\subfloat[NNN \label{fig:NNN911}]
		{\includegraphics[width = 1.5in]{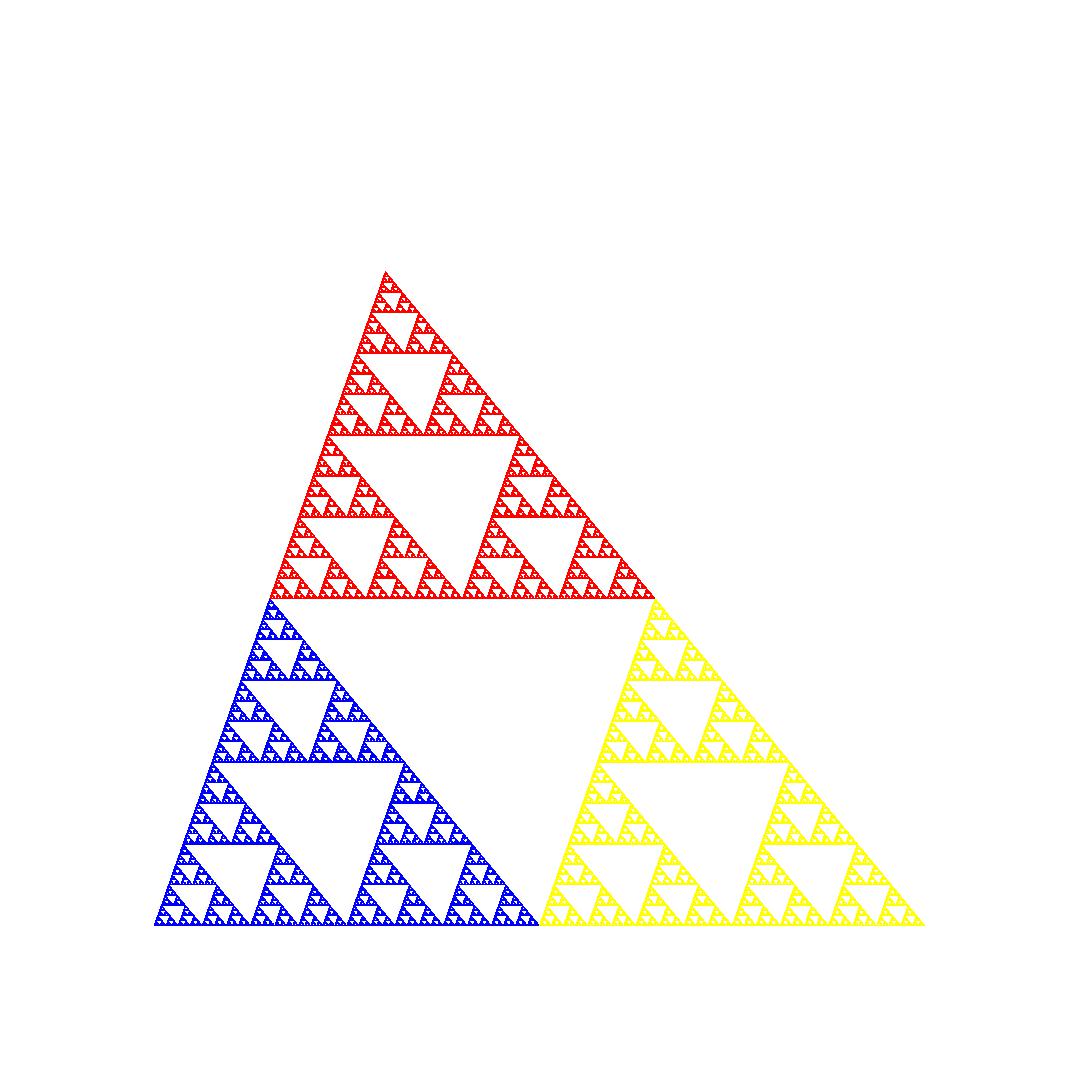}} &
		
		\subfloat[FNN \label{fig:FNN911}]
		{\includegraphics[width = 1.5in]{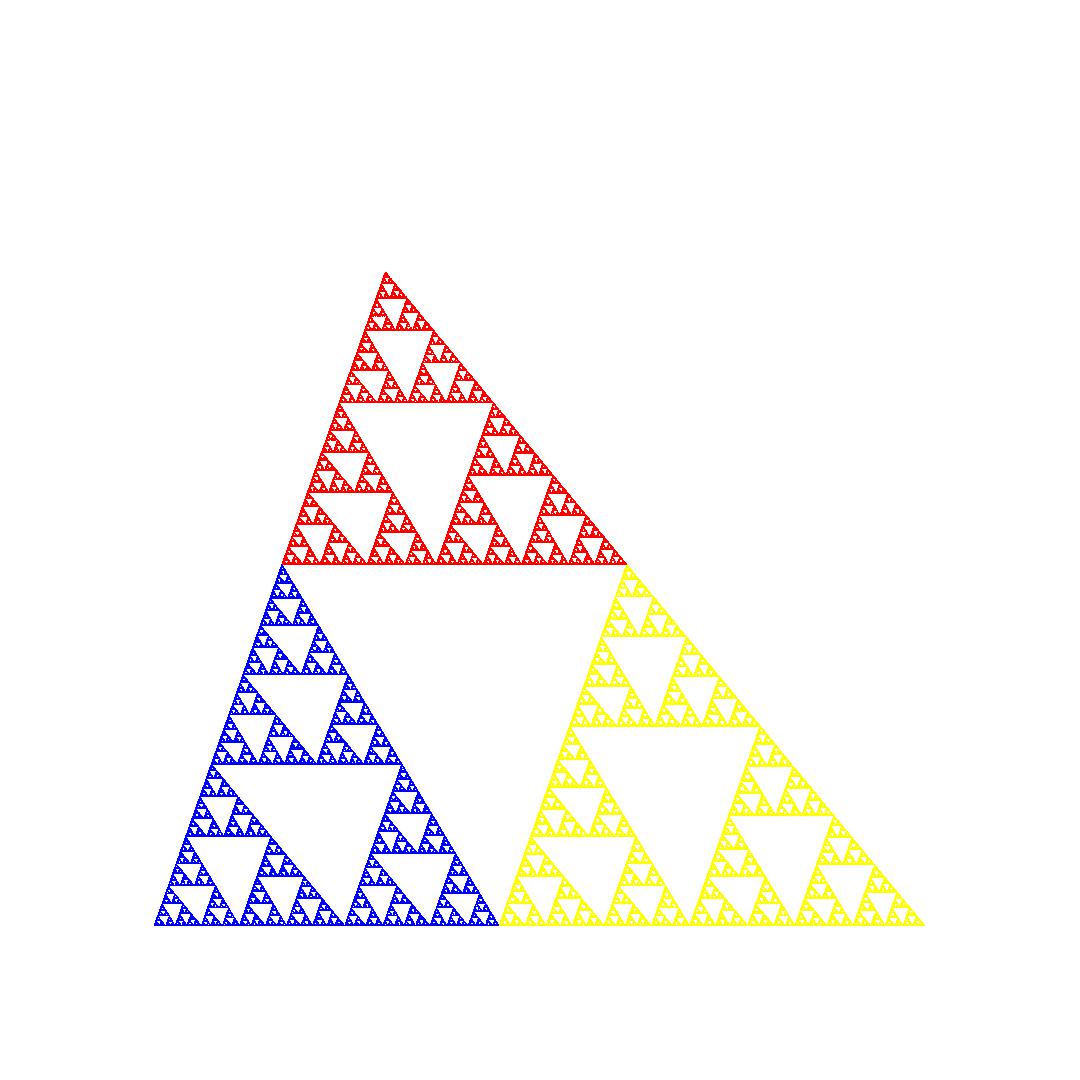}} &
	
		\subfloat[FFN \label{fig:FFN911}]
		{\includegraphics[width = 1.5in]{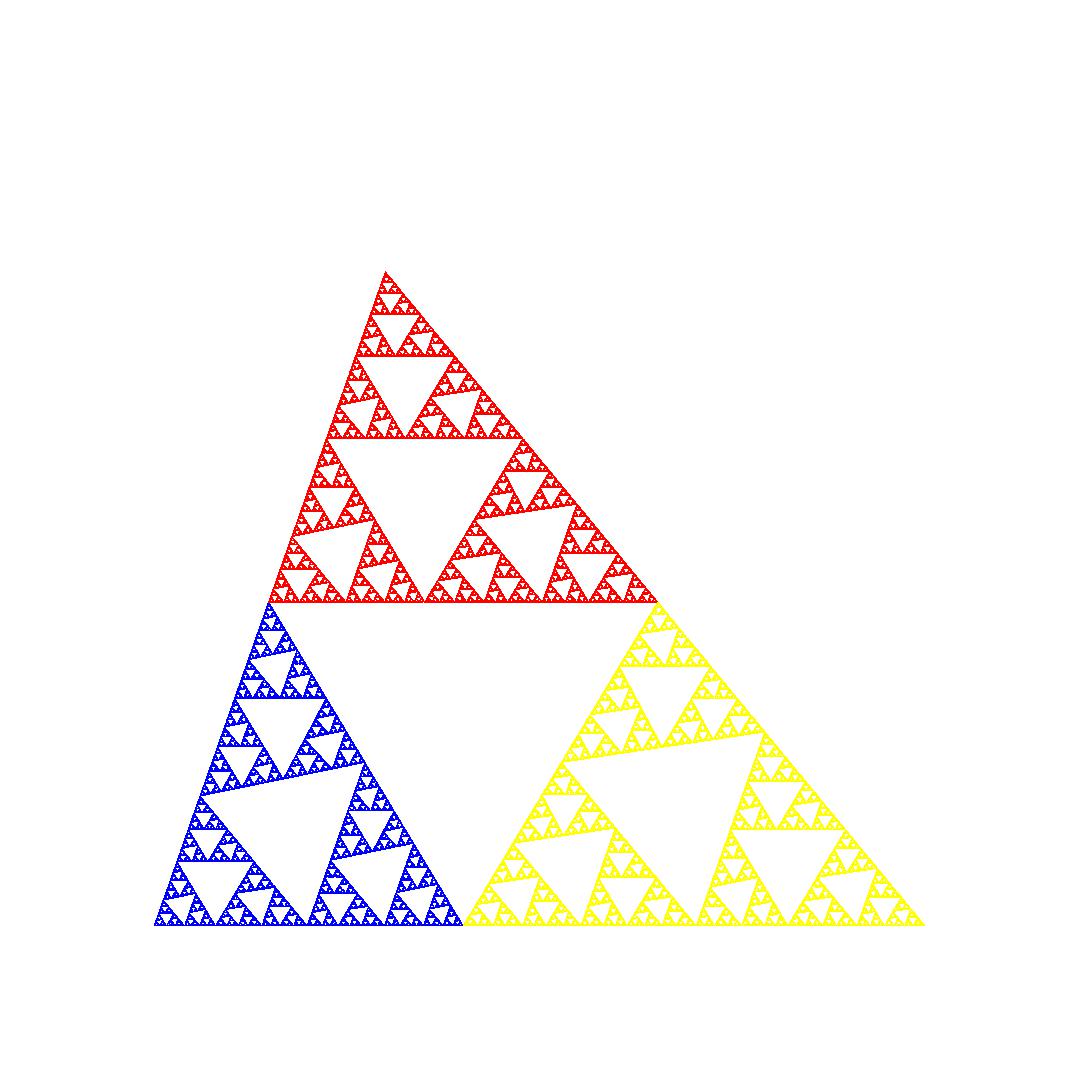}} &
		
		\subfloat[FFF \label{fig:FFF911}]
		{\includegraphics[width = 1.5in]{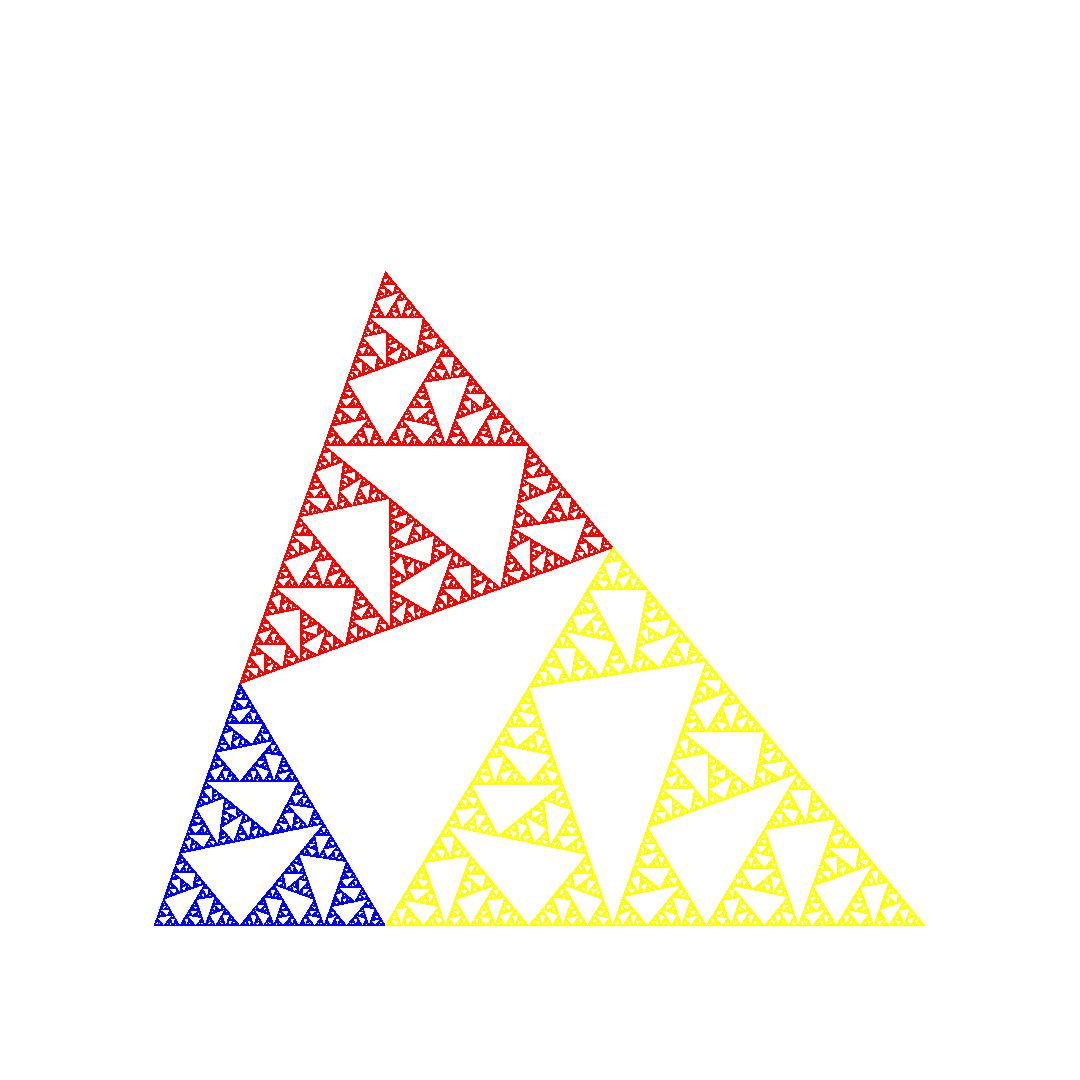}} 
	\end{tabular}
	\caption{Generalised Sierpinski triangles with a=1.1 and b=0.9}
	\label{fig:simplecon}
\end{figure}

\hspace*{\fill} \\
Let the side lengths of the fractal triangle be $a=1.1$ and $b=0.9$. By the similitudes defined in Section 3 the four types of generalised Sierpinski triangles are rendered in Figure \ref{fig:simplecon}. For each type of fractal triangle a different `hole' triangle is observed. It is seen that a non-flip map gives a side of the hole to be parallel to the exterior triangle while a flip map does not give a parallel line (in general). 
\vspace{-1cm}
It is seen in Section \ref{sec:mindimension} that the four types of generalised Sierpinski triangles do not, in general, have the same fractal dimension and this is seen in that the amount of region for where points exists is different for each triangle.
The fractal triangle can be thought of as the set $A$ such that $F(A) = \bigcup_{i=1}^N f_i(A) = A$ or it can be thought of as a filled in $\triangle ABC$ with a `hole' removed (the hole was previously described by $\triangle NMO$) with this process repeated ad infinitum on the resultant triangles. 

\vspace{2mm}
\subsection{Obtuse Generalised Sierpinski Triangles} \label{obtusesection} \hspace*{\fill} \\
A consequence of how the triangles were constructed is the extension of the maps that contain non-flips to produce obtuse generalised Sierpinski triangles . This is trivial for the classic Sierpinski triangles ($\triangle NNN$). However, the fractal Pedal ($\triangle FFF$) triangles can only be acute and right-angled; this can be seen in the fact that the they are constructed by removing the inner hole-triangle from the large triangle at each iteration. For obtuse triangles, the hole partially lies outside the initial obtuse triangle so this construction becomes nonsensical since it would have an attractor that does not resemble a triangle. 

\begin{figure}[H]
	\centering
	\begin{tabular}{ccc}
		\subfloat[NNN \label{fig:NNN142}]
		{\includegraphics[width = 1.8in]{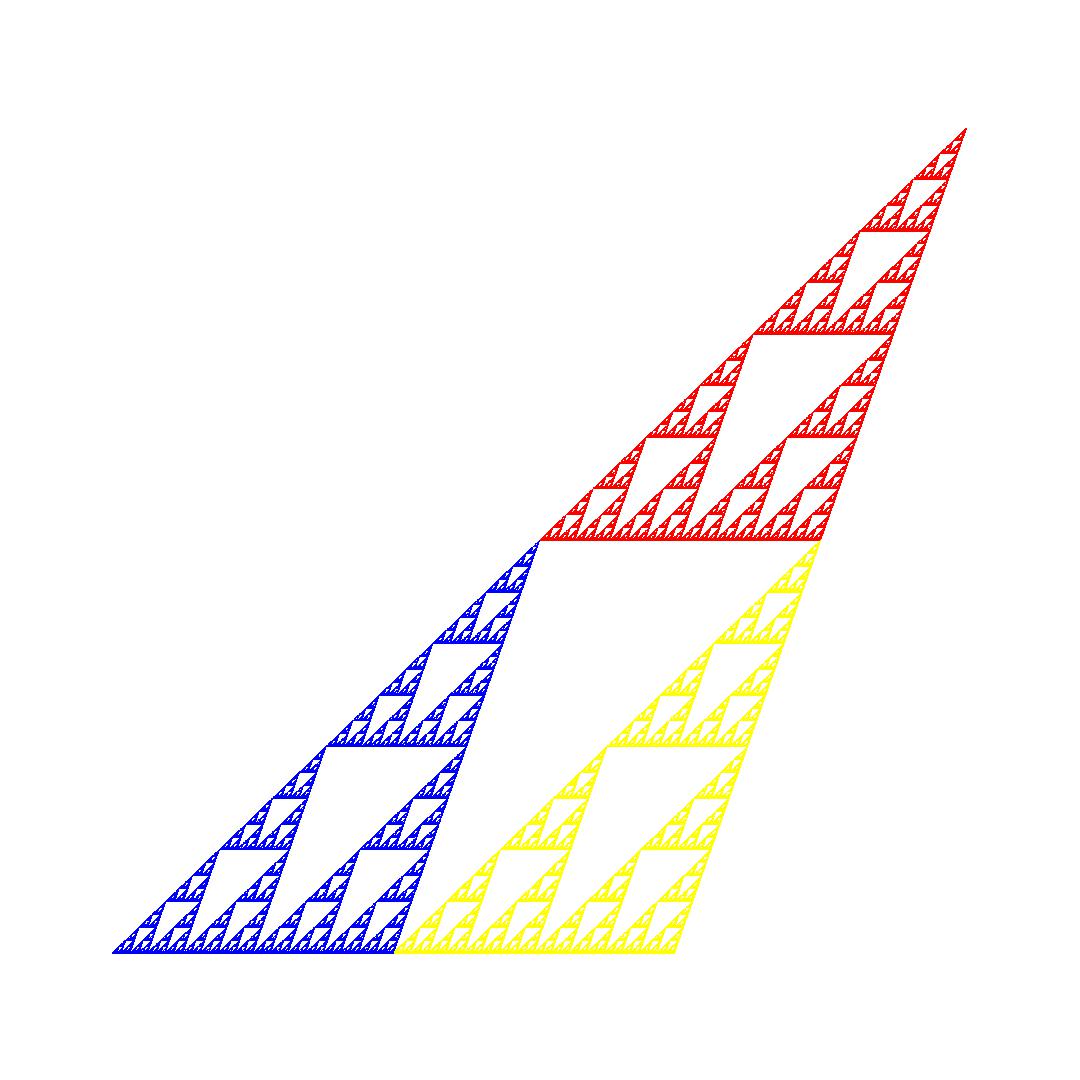}} &
		
		\subfloat[FNN \label{fig:FNN142}]
		{\includegraphics[width = 1.8in]{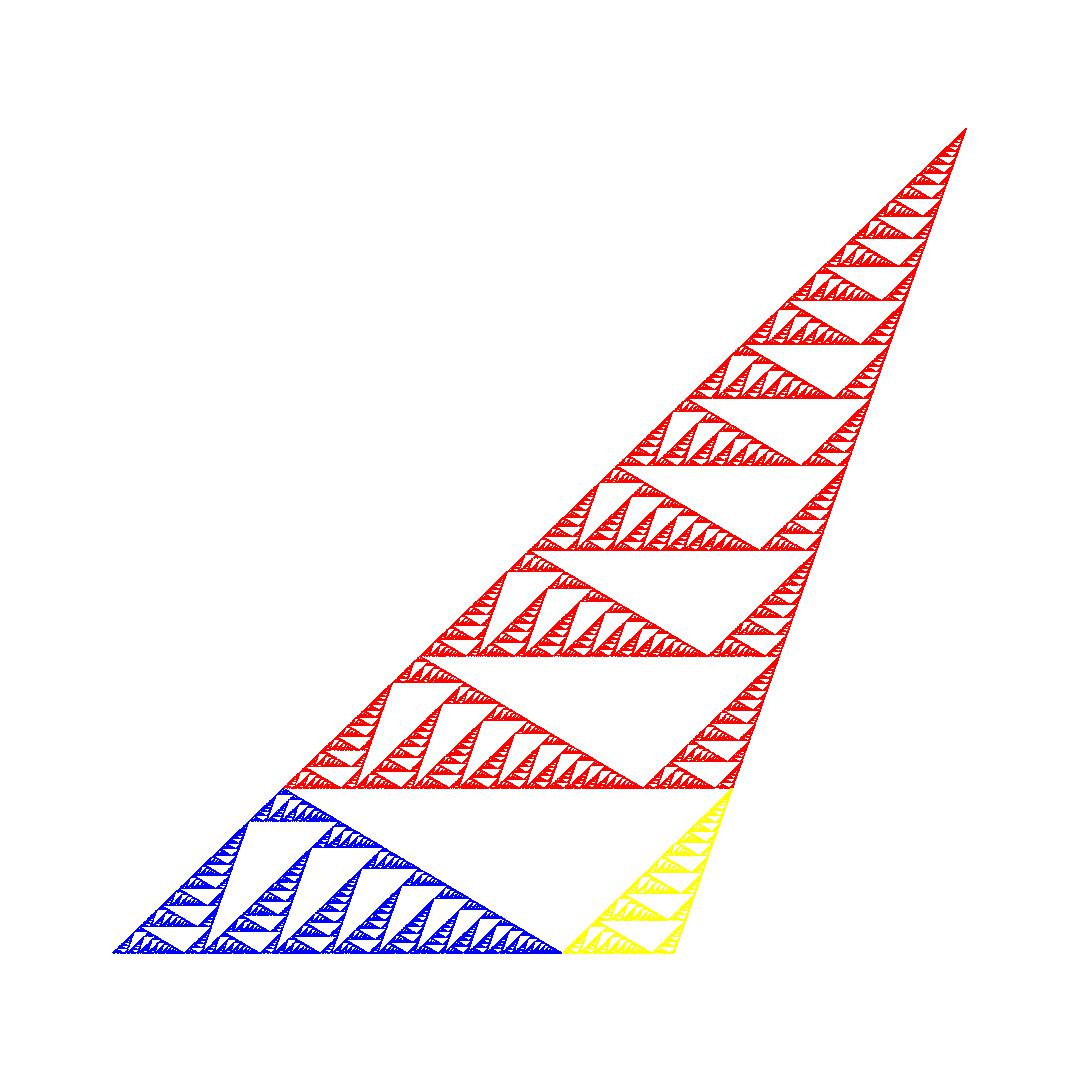}} &
		
		\subfloat[FFN \label{fig:FFN142}]
		{\includegraphics[width = 1.8in]{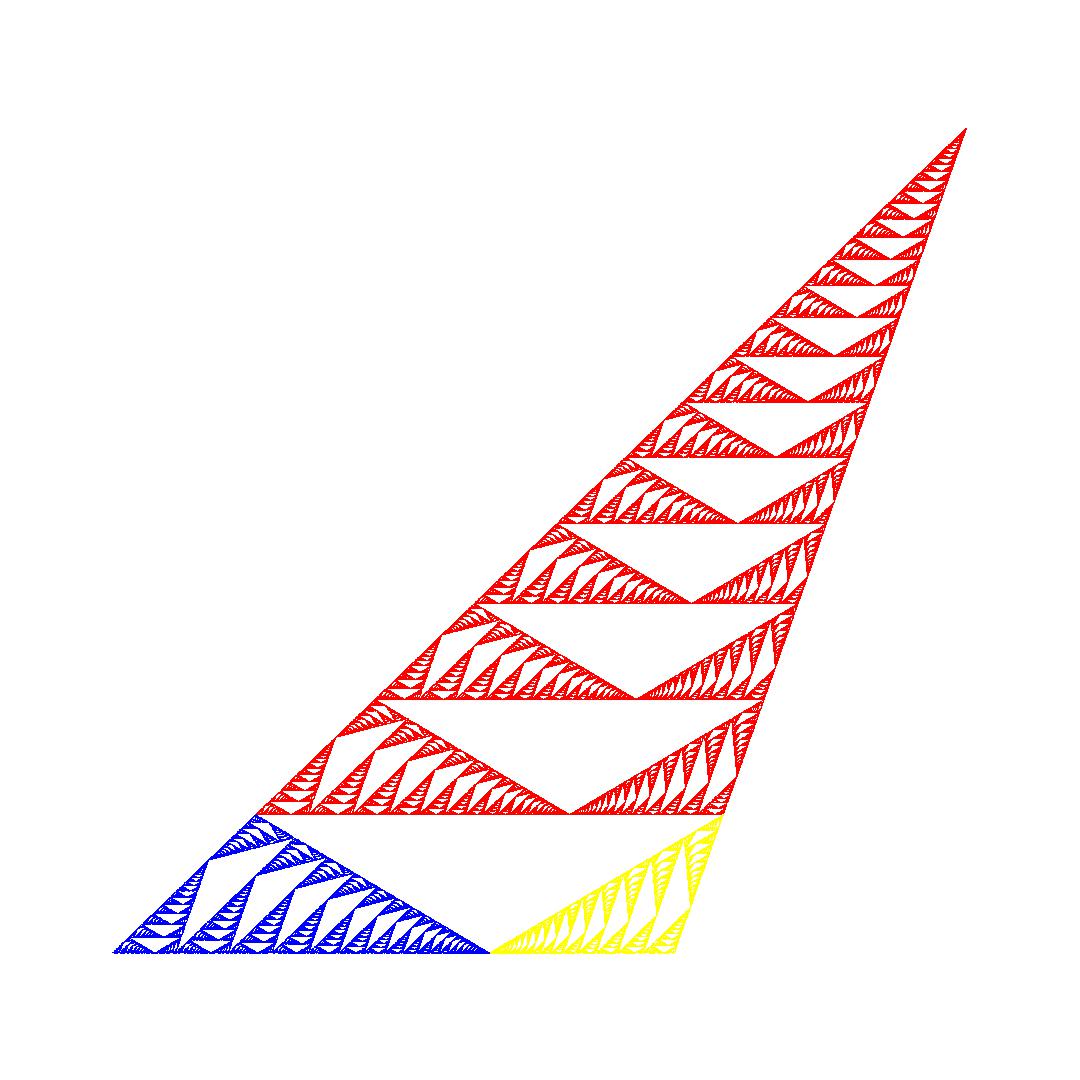}} 
	\end{tabular}
	\caption{Obtuse Generalised Sierpinski with a=1.4 and b=2}
\end{figure}
This extension to obtuse triangles also extends the range of the dimension of these triangles with the previously known lowest-dimension fractal triangle being the Sierpinski which has fixed dimension. Now, however, the $\triangle FFN$-type triangle is able to produce generalised Sierpinski triangles with dimension lower than Sierpinki triangle's dimension of $\frac{log(3)}{log(2)} \approx 1.58$.
\begin{figure}[H]
	\centering
	\includegraphics[width = 4.5in]{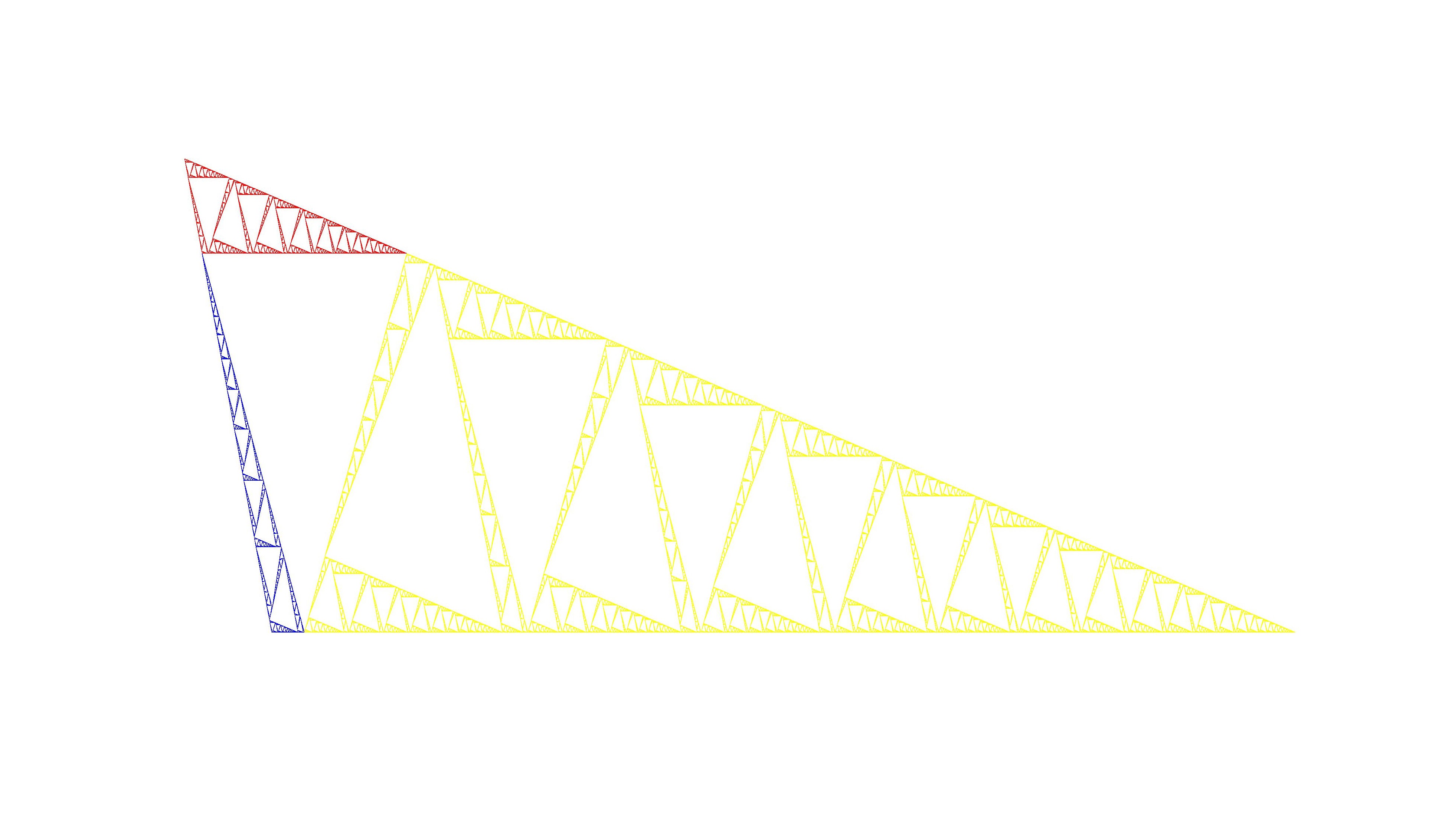}
	\caption{FFN triangle with $a=0.2, b=1.1$ and dimension $\approx 1.44$}
	\label{fig:FFN144}
\end{figure}

\section{\bf Minimising dimension of Generalised Sierpinski Triangles} \label{tridim} \label{sec:mindimension}
\subsection{Dimension of Sierpinski (NNN) Triangles}\hspace*{\fill} \\
The Moran--Hutchinson Theorem [\ref{hutmor}] can be used to calculate the fractal dimension $d$ of all Sierpinski triangles using the scaling ratios of the similitudes (which is $\alpha= \beta= \gamma= \frac{1}{2}$). The fractal dimension is the unique positive solution $d$ to:
\begin{equation*}
\alpha^d+\beta^d+ \gamma^d = 1 \quad \Rightarrow \quad 3\Big( \frac{1}{2} \Big) ^d = 1 \quad \Rightarrow \quad d = \frac{ln(3)}{ln(2)} \approx 1.58
\end{equation*}
Hence the dimension of the Sierpinski triangle is fixed at the known value of $\frac{log(3)}{log(2)}$.

\vspace{2mm}
\subsection{\bf Dimension of $\bf \triangle FNN$}\hspace*{\fill} \\
Firstly note that the scaling factors of our three maps in this situation are purely dependant on the $b$ side of the triangle. From Section \ref{fractrisec},

\begin{equation*}
\alpha = \frac{b}{b^{2}+1} \quad \beta = \frac{1}{b^{2}+1}  \quad \gamma=\frac{b^{2}}{b^{2}+1} 
\end{equation*} 
Therefore the system that needs to be minimised is:
\begin{equation}
\Big(\frac{b}{b^{2}+1}\Big)^d+\Big(\frac{1}{b^{2}+1}\Big)^d+ \Big(\frac{b^{2}}{b^{2}+1}\Big)^d = 1
\label{NNFdim}
\end{equation}
It will be shown that the dimension $d$ has a minimum at $b=1$. By implicit differentiation of Eq. \ref{NNFdim}, 
\begin{equation}
\frac{\mathrm{d}d}{\mathrm{d}b} =\frac{d \left(-\left(\frac{b}{b^2+1}\right)^d-2 \left(\frac{b^2}{b^2+1}\right)^d+b^2 \left(2
	\left(\frac{1}{b^2+1}\right)^d+\left(\frac{b}{b^2+1}\right)^d\right)\right)}{\left(b^3+b\right) \left(\left(\frac{1}{b^2+1}\right)^d \left(\log
	\left(\frac{1}{b^2+1}\right)\right)+\left(\frac{b}{b^2+1}\right)^d \log \left(\frac{b}{b^2+1}\right)+\left(\frac{b^2}{b^2+1}\right)^d \log
	\left(\frac{b^2}{b^2+1}\right)\right)}
\label{dbNNF}
\end{equation}
Substituting $b=1$ makes $\frac{\mathrm{d}d}{d\mathrm{d}b}=0$, which grants $b=1$ as a critical point of the function. It is clear for positive $d$ and $b$ that the denominator of Eq. \ref{dbNNF} is negative. To show $b=1$ is the global minimum it suffices to show the numerator is positive for $0<b<1$ and negative for $1<b<\infty $. 

$$\mathrm{numerator} = d \left(-\left(\frac{b}{b^2+1}\right)^d-2 \left(\frac{b^2}{b^2+1}\right)^d+b^2 \left(2 \left(\frac{1}{b^2+1}\right)^d+\left(\frac{b}{b^2+1}\right)^d\right)\right)$$
This is positive/negative if and only if
$$-b^d-2 \left(b^2\right)^d+b^2 \left(b^d+2\right)=2 \left(b^2-\left(b^2\right)^d\right)+\left(b^2-1\right) b^d$$
is positive/negative. We note that if $0<b<1$, $1< d\leq 2$, this equation  is positive and if $1<b< \infty$ and $1< d\leq 2$ it is negative then it would have been shown that $d$ was a global minimum.  If this were the case, $b=1$ is when the dimension obtains a minimum value, which makes $$\alpha=\beta=\gamma=\frac{1}{2} \quad \Rightarrow \quad d=\frac{log(3)}{log(2)}$$
Plotting Equation \ref{NNFdim} presents as a plausible conjecture that $d=\frac{log(3)}{log(2)}$ is the global dimension of this type of fractal triangle.
\begin{figure}[!ht]
	\centering
	\begin{tabular}{cc}
		\subfloat[Dimension of $\triangle FNN$ triangle vs. b ]
		{\includegraphics[width = 2in]{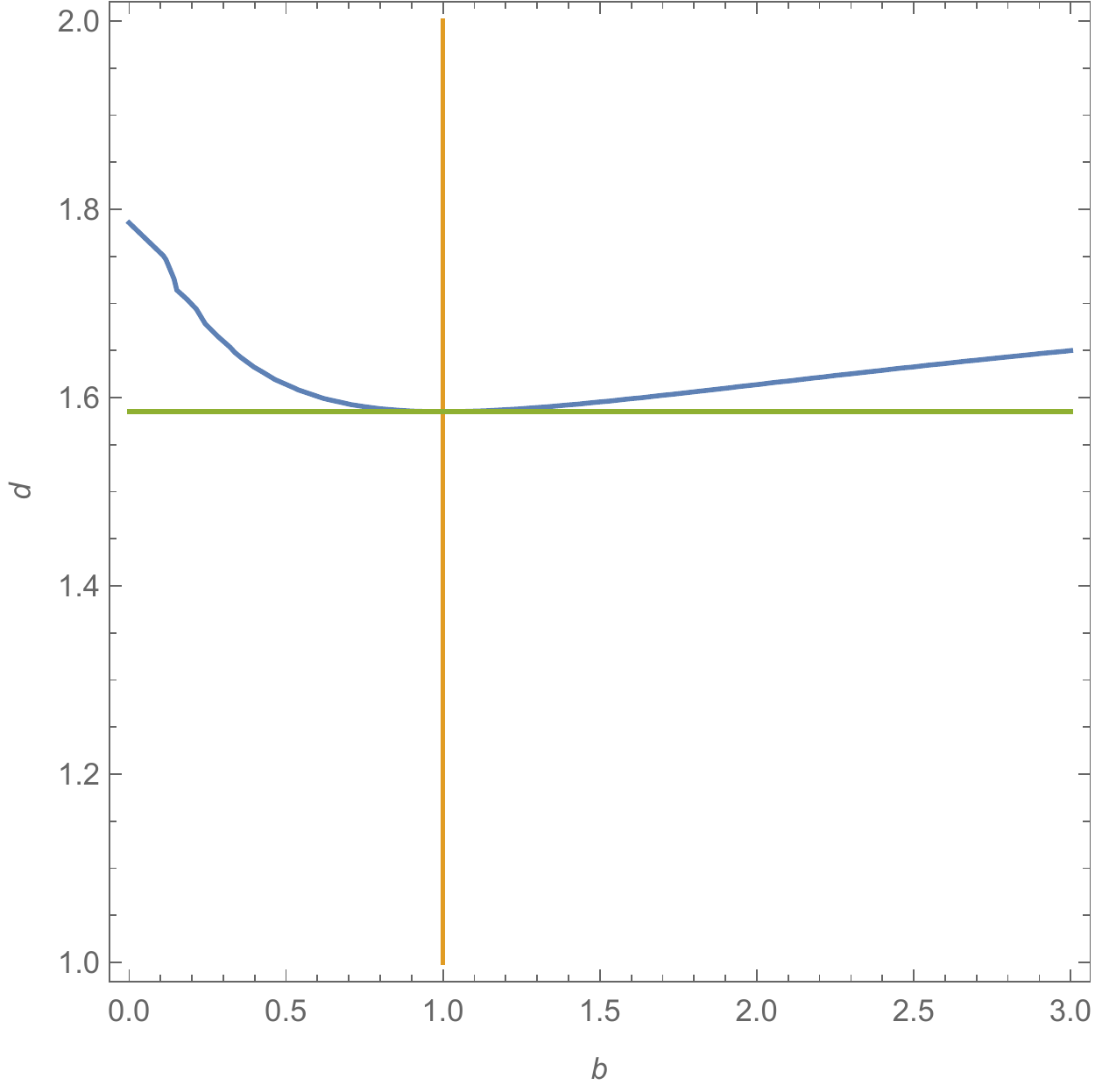}} \hspace{2cm}
		&
		\subfloat[Fractal dimension for $\triangle FFN$ against a and b \label{fig:FFNdim}]
		{\includegraphics[width = 2in]{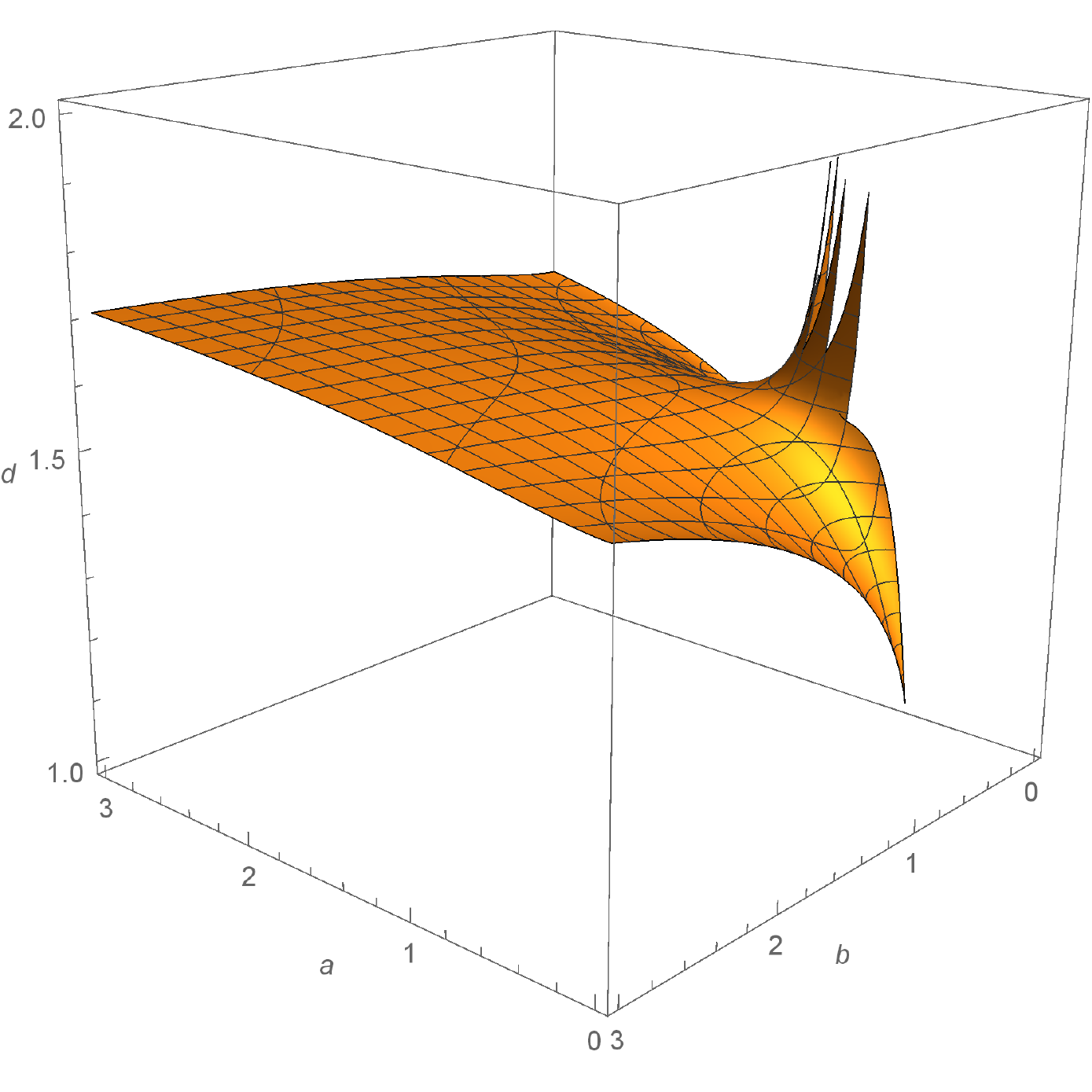}}
	\end{tabular}
	\caption{Dimension Plots for FFN \& FNN}
\end{figure}

\vspace{2mm}
\subsection{\bf Dimension of $\bf \triangle FFN$}\hspace*{\fill} \\
From Section \ref{fractrisec}, the three scaling factors for this IFS are: 
\begin{equation*}
\alpha = \frac{b}{a^2+b^{2}} \quad \beta = \frac{a}{a^2+b^{2}}  \quad \gamma=\frac{a^2+b^{2}-1}{a^2+b^{2}}
\end{equation*} 
Thus the equation for the fractal dimension is:

\begin{equation*}
\Big( \frac{b}{a^2+b^{2}}\Big)^d +  \Big(\frac{a}{a^2+b^{2}}\Big)^d  +\Big(\frac{a^2+b^{2}-1}{a^2+b^{2}}\Big)^d=1
\label{FFNdim}
\end{equation*} 
Using the method of Lagrange multipliers, the critical points of this function are tested. First we define the function $\Lambda (a,b,d,\lambda)$:
\begin{equation}
	\Lambda (a,b,d,\lambda) := d - 	\lambda \Big(\Big( \frac{b}{a^2+b^{2}}\Big)^d +  \Big(\frac{a}{a^2+b^{2}}\Big)^d  +\Big(\frac{a^2+b^{2}-1}{a^2+b^{2}}\Big)^d -1\Big)
\end{equation}
Taking the partial derivative with respect to each of the variables yields,
\begin{align*}
    \Lambda_a &= -\frac{\lambda \left(d (b-a) (a+b) \left(\frac{a}{a^2+b^2}\right)^{d-1}-2 a b d \left(\frac{b}{a^2+b^2}\right)^{d-1}+2 a d
    	\left(1-\frac{1}{a^2+b^2}\right)^{d-1}\right)}{\left(a^2+b^2\right)^2} \\   
	\Lambda_b&=-\frac{\lambda \left(-2 a b d \left(\frac{a}{a^2+b^2}\right)^{d-1}+d (a-b) (a+b) \left(\frac{b}{a^2+b^2}\right)^{d-1}+2 b d
		\left(1-\frac{1}{a^2+b^2}\right)^{d-1}\right)}{\left(a^2+b^2\right)^2}\\	
	\Lambda_d&=1-\lambda \Big(\left(\frac{a}{a^2+b^2}\right)^d \log \left(\frac{a}{a^2+b^2}\right)+ \left(\frac{b}{a^2+b^2}\right)^d \log
	\left(\frac{b}{a^2+b^2}\right)\\
	&+\left(1-\frac{1}{a^2+b^2}\right)^d \log \left(1-\frac{1}{a^2+b^2}\right)\Big) \\	
	\Lambda_\lambda&=-\left(\frac{a}{a^2+b^2}\right)^d-\left(\frac{b}{a^2+b^2}\right)^d-\left(1-\frac{1}{a^2+b^2}\right)^d+1
\end{align*}
Substituting $a=b=1$ and choosing $\lambda=\frac{-1}{\log(2)}$ into all of these makes them vanish, which grants $(a,b,\lambda) = (1,1,\frac{-1}{\log(2)})$ as a critical point of the function. When solving for $d$, this critical point yields $d=\frac{log(3)}{log(2)}$. This value is seen to be a critical point of all the other types of generalised Sierpinski triangles. But we show that this critical point is not a minimum: as noted above in Figure \ref{fig:FFN144}, there are other points that have dimension lower than this value. Plotting Equation \ref{FFNdim} shows this critical point is actually a saddle point.

\begin{figure*}[!ht]
 	\centering
 	\includegraphics[width=3in]{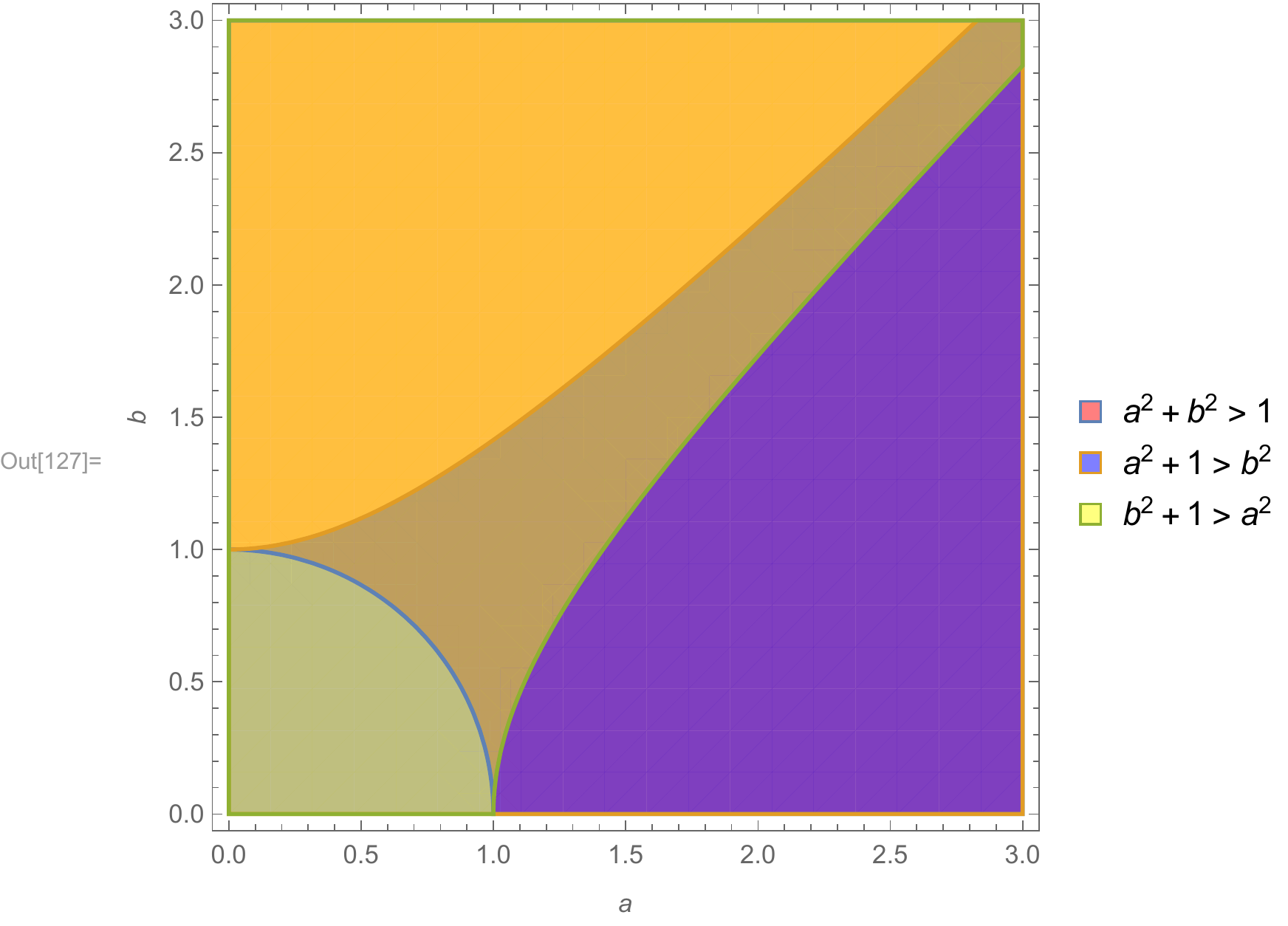}
 	\caption{Regions in a and b for different types of triangle}
 	\label{acreg}
 \end{figure*}

By inspection of Figure \ref{FFNdim}, the function is symmetric along the plane $a=b$ so it suffices to analyse only one side of this plane. The function tends to zero as $b \rightarrow 1$ and $a \rightarrow 0$ when approaching from the `orange' region of Equation \ref{acreg}, which makes sense geometrically as this would mean the triangle flattening to a straight line. But if the point $b=1 ,a=0$ was approached along the curve $a^2+b^2=1$ (geometrically, this is the sequence of right-angle triangles with hypotenuse $1$), the function would equal 2. This can be seen by solving the equation for dimension along this curve.
\begin{align*}
	\Big( \frac{b}{a^2+b^{2}}\Big)^d +  \Big(\frac{a}{a^2+b^{2}}\Big)^d  +\Big(\frac{a^2+b^{2}-1}{a^2+b^{2}}\Big)^d&=1 ,\quad a^2+b^2=1\\
	\Big( \frac{b}{1}\Big)^d +  \Big(\frac{a}{1}\Big)^d  +\Big(\frac{1-1}{1}\Big)^d&=1\\
	a^d+b^d&=1\\
	\Rightarrow d&=2 \quad \text{as} \quad a^2+b^2=1
\end{align*}
But at the point $a=0$ and $b=1$, the equation gives $1^d=1$ so $d$ cannot be solved in this way. The limits approaching this point are not equal from all directions and hence $(a,b) = (0,1)$ is a discontinuity. Geometrically, if the point is approached with an obtuse triangle (i.e the orange region of Figure \ref{acreg}) the limit of $d$ is one, and if the point is approached with a right-angle or acute triangle (the brown region of Figure \ref{acreg}) the limit of $d$ is two at $a=0, b=1$.
 \\
Therefore, obtuse triangles can produce generalised Sierpinski triangles of dimension arbitrarily close to one and the acute or right-angle triangles may produce generalised Sierpinski triangles with lesser dimension than Sierpinski's dimension of $\frac{log(3)}{log(2)}$. From standard calculus minimisation techniques, a continuous function (which the dimension equation os when $\triangle FFN$ is restricted to the acute region) that is also bounded attains a maximum/minimum on the boundary or at its critical points. Locally in the plotted region, the function appears to not have any other critical points other than the saddle point $a=b=1$ so any maximums/minimums may lay on the boundary. It has already been established that $d=2$ along the boundary $a^2+b^2=1$ and is a global maximum as the space is $\mathbb{R}^2$. In terms of the minimums for the bounded local acute area, along the boundary $a^2=b^2+1$, the function goes below the Sierpinski's dimension. 
\vspace{2mm}
\subsection{\bf Dimension of Pedal ($\bf \triangle FFF$) Triangles} \hspace*{\fill} \\
The Pedal triangles have already been well investigated in sources \cite{PedalCastellano}, \cite{PedalDing}, \cite{PedalDingNC} and \cite{PedalHitt} . It has been established that for acute triangles the minimum dimension is the Sierpinski triangle dimension of $\frac{log(3)}{log(2)}$. By the process of Lagrange multipliers and as demonstrated for $\triangle FFN$ triangles with the explicitly defined maps, it can be shown that $a=b=1$ is a critical point of the dimension function and that these side lengths give the Sierpinski dimension for the Pedal triangle. As there are no sensible constructions for obtuse fractal pedal triangles, it may be concluded that the global minimum of the fractal dimension for the Pedal ($\triangle FFF$) triangle is the Sierpinski dimension.

\vspace{6mm}
\section{\bf Tiling of $\bf \mathbb{R}^2$ by Generalised Sierpinski Triangles} \label{fractalblowups}
The purpose of this section is to use the newly discovered triangles to create fractal tilings of $\mathbb{R}^2$ with fewer symmetries than those tilings made from the standard Sierpinski triangle. As each of the generalised Sierpinski triangle IFSs are similitudes with attractors that obey the open set condition (Defn \ref{OSC}), their dimensions are well-defined and can be used to explicitly investigate the fractal tilings they create from the results found in \cite{BarnNew}.\\

\vspace{2mm}
\subsection{\bf Tiling Notation} \hspace*{\fill} \\
The following notation is taken from \cite{BarnTile}. \\
A tile is a perfect compact non-empty subset of $\mathbb{R}^2$. It is immediate that all of the generalised Sierpinkski triangles are themselves tiles. A tiling of $\mathbb{R}^2$ is a union of tiles all with equal Hausdorff dimension such that no tiles are overlapping. The set of all tilings denoted $\mathbb{T}$. To construct a tiling from a generalised Sierpinski triangle, the following notation is introduced:\\
Let $[3]=\{1,2,3\}$ and $[3]^k$ be the set of words of length $k$ from the alphabet $[3]$. Then denote the set of all finite words $[3]^{*}=\cup_{k=0}^{\infty}[3]^k$, where $[3]^0$ is the empty word and the set of infinite words is $[3]^{\infty}$ .\\
For $\omega \in [3]^{\infty}$ and $k \in \mathbb{N}$ the following notation is used,
$$ \omega | k := \omega_{1}\omega_{2}\omega_{3}...\omega_{k}$$
Where $\omega_{i}\in [3]$. For an IFS  $F=\{\mathbb{R}^2 ; f_1, f_2, f_3\}$ with attractor $A$, then $f_{\omega | k}$ is defined,
$$ f_{\omega | k} := f_{\omega_{1}} \circ f_{\omega_{2}} \circ f_{\omega_{3}} \circ ...f_{\omega_{k}} (A)$$ 
and furthermore as all the generalised Sierpinski triangles are invertible,
$$ (f^{-1})_{\omega | k}:=  f^{-1}_{\omega_{1}} \circ  f^{-1}_{\omega_{2}} \circ  f^{-1}_{\omega_{3}} \circ ... f^{-1}_{\omega_{k}} (A)$$ 
Now for $\theta \in [3]^{\infty}$ and $\omega \in [3]^{k}$, a tile can be defined as,
$$t_{\theta,\omega}:=((f^{-1})_{\theta|k}\circ f_{\omega})(A)$$
and a tiling being the union of created tiles,
$$T_{\theta,k} : = \{t_{\theta,\omega} : \omega \in[3]^{k} \}$$
From \cite{BarnTile}, taking the union over all $k$ yields an unbounded tiling denoted
$$T_{\theta} := \bigcup_{k=1}^{\infty} T_{\theta,k}$$
For each such tiling, a prototile set  $\mathcal{T}$ is a minimal set of tiles such that every tile in $T_{\theta}$ is an isometric copy of a tile in $\mathcal{T}$.

\vspace{2mm}
\subsection{\bf Periodic Tilings and Finite Prototile Sets}

The work completed in \cite{BarnNew} shows that if an IFS obeys the  algebraic condition that the scaling ratios are all positive integer powers of a positive number $s$, then the tilings constructed from that IFS will be quasiperiodic and have a prototile set which contains finitely many tiles. For generalised Sierpinski triangles with scaling ratios $\alpha , \beta , \gamma$, the algebraic condition is explicitly $\alpha=s^{a_1} , \beta=s^{a_2} , \gamma=s^{a_3}$ and furthermore the prototile set is $\mathcal{T} = \{sA,s^2A,...,s^{a_{max}}A\}$ where $a_{\max}=\max\{a_1,a_2,a_3\}$. This is equivalent to requiring $\gamma = \alpha^x = \beta^y$ where $x$ and $y$ are integers and the selection $\alpha^{\frac{1}{y}}$ can be made to assure the condition.  As there are explicit formulas for the scaling factors of  $\alpha , \beta , \gamma$, examples of periodic tilings can be directly constructed.

\vspace{2mm}

\bf{$\triangle FFF$ Example}
\hspace*{\fill} \\
\mdseries
Using the above formulation and setting $x=y=2$ in the case of the $\triangle FFF$ triangle yields $s=\alpha$ and $\mathcal{T} = \{\alpha A, \alpha^2 A\}$. Solving this set of equations for $b$ and $a$ is easily done by:

	\begin{align*}
		\alpha^2&=\beta^2 \Rightarrow a=b\\
		\Rightarrow \left(\frac{1}{2a}\right)^2&=\frac{2a^2-1}{2a^2}\\
		\Rightarrow 1&=4a^2-2\\
		\Rightarrow a&=b=\frac{\sqrt{3}}{2}
	\end{align*}
	
Therefore the isosceles $\triangle FFF$ triangle with the above side lengths grants a periodic tiling out of two tile sizes isomorphic to the tiles in $\mathcal{T}$. Creating this tiling with $\theta = 12121212...\in [3]^{\infty}$ zoomed in at the point $(0.5,0)$ of $\mathbb{R}^2$ yields Figure \ref{FFFtile}:
\begin{figure*}[!ht]
	\centering
	\includegraphics[width=6.5in]{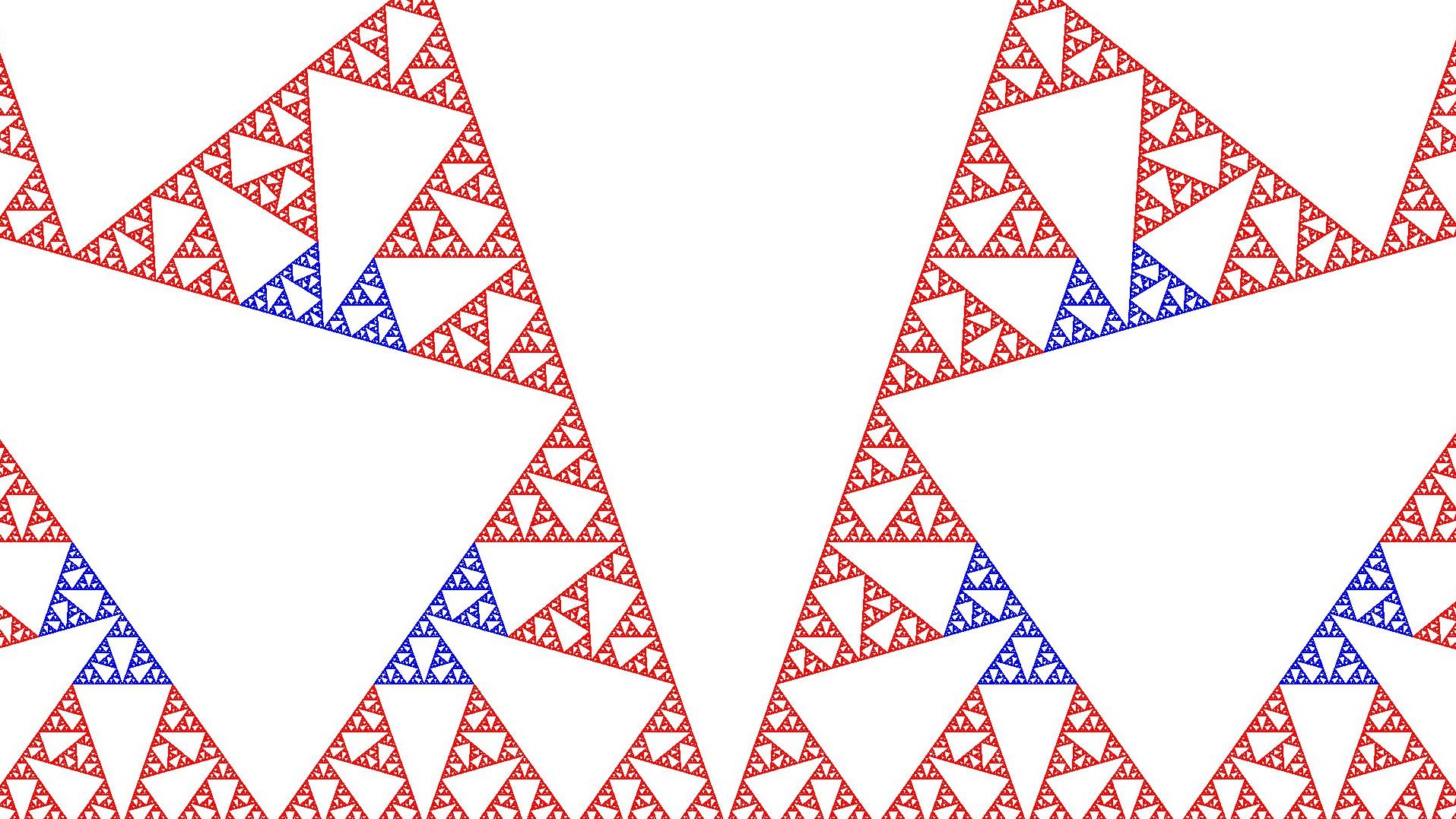}
	\caption{$T_{12121212...}$ at (0.5,0)}
	\label{FFFtile}
\end{figure*}\\
where Fig:~\ref{FFFtile} has the isomorphic tiles of the prototile set coloured in red ($\alpha A$) and blue ($\alpha^2 A$).\\
\vspace{2mm}

\pagebreak

\end{document}